\documentclass[11pt]{article}
\usepackage{graphicx,amssymb,amsmath,algorithmic}
\usepackage[all]{xy}
\newtheorem{theorem}{Theorem}

\newtheorem{lemma}[theorem]{Lemma}

\newtheorem{corollary}[theorem]{Corollary}

\def\EE{{\mathbb E}}
\def\M{{\mathbb M}}
\def\PP{{\mathbb P}}

\def\deg{{\mathrm{deg}}}

\def\ind{{\mathbf{1}}}
\newcommand\bp[1]{\noindent {\em Proof{#1}.} $\quad$}
\def\ep{\hfill $\Box$}

\def\R{{\mathbb R}}
\def\C{{\mathbb C}}
\def\N{{\mathbb N}}

\def\cG{{\mathcal G}}

\def\cS{{\mathcal S}}
\def\HH{{\mathbb H}}

\def\cM{{\mathcal M}}
\def\cH{{\mathcal H}}

\def\cU{{\mathcal U}}

\def\cP{{\mathcal P}}

\def\cN{\mathcal{N}}
\def\hpi{\widehat{\pi}}
\def\hphi{\widehat{\phi}}
\def\cGr{{{\mathcal G}_d}}
\def\cR{{\mathcal R}}

\title{Matchings on infinite graphs}
\author{Charles Bordenave \footnote{Institut de Math\'ematiques -
    Univ. de Toulouse \& CNRS  - France. Email:
    bordenave@math.univ-toulouse.fr } \,  Marc
  Lelarge\footnote{INRIA-\'Ecole Normale Sup\'erieure - France. Email:
    marc.lelarge@ens.fr} \,and Justin Salez\footnote{INRIA-\'Ecole
    Normale Sup\'erieure - France. Email: justin.salez@ens.fr}}
\date{}

\begin{document}
\maketitle
\begin{abstract}
Elek and Lippner (2010) showed that the convergence of a sequence of bounded-degree graphs implies the existence of a limit for the proportion of vertices covered by a maximum matching. We provide a characterization of the limiting parameter via a local recursion defined directly on the limit of the graph sequence. Interestingly, the recursion may admit multiple solutions, implying non-trivial long-range dependencies between the covered vertices. We overcome this lack of correlation decay by introducing a perturbative parameter (temperature), which we let progressively go to zero. This allows us to uniquely identify the correct solution. In the important case where the graph limit is a unimodular Galton-Watson tree, the recursion simplifies into a distributional equation that can be solved explicitly, leading to a new asymptotic formula that considerably extends the well-known one by Karp and Sipser for Erd\H{o}s-R\'enyi random graphs. 
\end{abstract}


\section{Introduction}
A \textit{matching} on a finite graph $G=(V,E)$ is a subset of
pairwise non-adjacent edges $M\subseteq E$. The $|V|-2|M|$ isolated
vertices of $(V,M)$ are said to be \textit{exposed} by $M$. We let
$\M(G)$ denote the set of all possible matchings on $G$. The
\textit{matching number} of $G$ is defined as
\begin{equation}
\label{eq:mn}
\nu(G)=\max_{M\in\M(G)}|M|,
\end{equation}
and those $M$ which achieve this maximum -- or equivalently, have the
fewest exposed vertices -- are called \textit{maximum matchings}. The \textit{normalized matching number} of $G$ is simply $\nu(G)/|V|$. 

Our results belong to the theory of convergent graph sequences. Convergence of bounded degree graph sequences was defined
by Benjamini and Schramm \cite{bensch}, Aldous and Steele
\cite{aldste}, see also Aldous and Lyons \cite{aldlyo}. The notion of local weak convergence has
then inspired a lot of work \cite{BenjaminiSS08}, \cite{res10}, \cite{borgs2010left},
\cite{elek}, \cite{lyotree}. In \cite{eleklippner}, it is shown that the convergence of a sequence of bounded-degree graphs guarantees the existence of a limit for their normalized matching numbers. We provide a characterization of the limiting parameter via a local recursion defined directly on the limit of the graph sequence. In the important case where the graph limit is a unimodular Galton-Watson tree, the recursion simplifies into a distributional equation that can be solved explicitly, leading to a new asymptotic formula.

A classical example in this context is the Erd\H{o}s-R\'enyi random graph with average degree $c$ on $n$ vertices, denoted by $G(n,c/n)$ : as $n\to\infty$, $G(n,c/n)$ converges in the local weak sense to a Galton-Watson tree with degree distribution Poisson with parameter $c$. In this case, Karp and Sipser \cite{karpsipser} showed that
\begin{equation}
\label{eq:KS}
\frac{\nu(G(n,c/n))}{n}\xrightarrow[n\to\infty]{}1-\frac{t_c+e^{-ct_c}+ct_ce^{-ct_c}}{2},
\end{equation}
where $t_c\in(0,1)$ is the smallest root of $t=e^{-ce^{-ct}}$ (we will see in the sequel that the convergence is almost sure). 
The explicit formula (\ref{eq:KS}) rests on the analysis of a heuristic
algorithm now called Karp-Sipser algorithm.
The latter is based on the
following observation : if $e\in E$ is a pendant edge (i.e. an edge incident to a vertex of degree one) in $G=(V,E)$, then
there is always a maximum matching that contains $e$, so all edges that are
adjacent to $e$ may be deleted without affecting $\nu(G)$. The first stage
of the algorithm consists in iterating this until no more pendant edge is
present. This is the \textit{leaf-removal process}.
$G$ is thus simplified into a sub-graph with only isolated
vertices, matched pairs, and a so-called \textit{core} with minimum degree
at least $2$. As long as that core is non-empty, one of its edges is
selected uniformly at random, the adjacent edges are deleted, and the whole
process starts again. When the algorithm stops, the remaining edges clearly
form a matching on $G$, but its size may be far below $\nu(G)$ due to the
sub-optimal removals on the core.

On $G(n,c/n)$, the dynamics of the deletion process can be approximated in the $n\to \infty$ limit by differential equations which can be explicitly solved.
In particular, the asymptotic size of both the optimal part
constructed in the first stage, and the sub-optimal part constructed on the
core can be evaluated up to an $o(n)$ correcting term (which has been later
refined, see \cite{ksrevisited}). Moreover, the second part happens to be almost
perfect, in the sense that only $o(n)$ vertices are exposed in the
core. This guarantees that the overall construction is asymptotically
optimal, and the asymptotic formula for $\nu(G_n)$ follows.
More recently, the same technique has
been applied to another class of random graphs with a fixed log-concave degree profile \cite{friezebohman}, resulting in the asymptotical existence of an \textit{almost perfect matching} on these graphs :
\begin{equation}
\label{eq:BF}
\frac{\nu(G_{n})}{|V_n|}\xrightarrow[n\to\infty]{}\frac{1}{2}.
\end{equation}
In both cases, the proof of optimality -- and hence
the asymptotic formula for $\nu(G_n)$ -- relies on the fact that the second stage exposes only $o(n)$ vertices, which is bound to fail as soon as one considers more general graph ensembles where the core does
not necessarily admits an almost-perfect matching. We give
simple examples in the Appendix. By using a completely different approach -- namely establishing and solving an appropriate \textit{recursive distributional equation} (a usual ingredient of the objective method, see \cite{aldousbandyopadhya}) --, we manage to obtain a general formula that considerably extends the above results. 

The rest of our paper is organized as follows: we state our main results in Section \ref{sec:res}.
In Section \ref{sec:z}, we extend the Boltzmann-Gibbs distribution over matchings on a finite graph to infinite graphs. This will allow us to derive our Theorem \ref{th:main} in Section \ref{sec:zero}.
We deal with the specific cases of trees (and random graphs) in Section \ref{sec:rdet}. 
We end the paper with an Appendix presenting simple examples of graphs for which the limiting local recursion admit multiple solutions.

\section{Results}\label{sec:res}

Let us start with a brief recall on local weak convergence (see \cite{bensch, aldste} for details). A \textit{rooted graph} $(G,\circ)$ is a graph $G=(V,E)$ together with the
specification of a particular vertex $\circ\in V$, called the
\textit{root}. We let $\cG$ denote the set of all locally finite connected rooted graphs considered up to \textit{rooted isomorphism}, i.e. $(G_1,\circ_1)\equiv(G_2,\circ_2)$ if there exists a bijection $\gamma\colon V_1\to V_2$ that preserves roots ($\gamma(\circ_1)=\circ_2$) and 
adjacency ($uv\in E_1\Longleftrightarrow \gamma(u)\gamma(v)\in E_2$). In the space $\cG$, a sequence $\left\{(G_n,\circ_n);n\in \N\right\}$ \textit{converges locally} to $(G,\circ)$ if for every radius $k\in\N$, there is $n_k\in\N$ such that
$$n\geq n_k\Longrightarrow [G_n,\circ_n]_k\equiv [G,\circ]_k.$$ 
Here, $[G,\circ]_k$ denotes the finite rooted subgraph induced by the vertices lying at graph-distance at
most $k$ from $\circ$. It is not hard to construct a distance which metrizes this notion of convergence and turns $\cG$ into a complete separable metric space. We can thus import the usual machinery of weak convergence of probability measures on Polish spaces (see e.g. \cite{billingsley}). We define $\cP(G)$ as the set of probability measure on $\cG$. There is a natural procedure for turning a finite deterministic
graph $G=(V,E)$ into a random element of $\cG$ : one simply chooses uniformly at random a vertex
$\circ\in V$ to be the root, and then restrains $G$ to the
connected component of $\circ$. The resulting law is denoted by $\cU(G) \in \cP(\cG)$. If $(G_n)_{n\in\N}$ is a sequence of finite graphs
such that $(\cU(G_n))_{n\in\N}$ admits a weak limit $\rho\in\cP(\cG)$, we call $\rho$ the \textit{random weak limit} of the sequence  $(G_n)_{n\in\N}$. Finally, for any $d\geq 0$, we define $\cGr$ as the space of all rooted connected graphs with maximal degree no more than $d$.

Rather than just graphs $G=(V,E)$, it will be sometimes convenient to
work with \textit{discrete networks} $G=(V,E,\cM)$, in which the additional
specification of a \textit{mark map} $\cM\colon E\to\N$ allows to attach useful local
information to edges, such as their absence/presence in a certain matching. We then simply require the isomorphisms in the
above definition to preserve these marks.

The first main implication of our work is that the local weak convergence of a sequence of graphs is enough to guarantee the convergence of their normalized matching numbers to a quantity that can be described directly on the random weak limit of the graph sequence. 
\begin{theorem}
\label{th:main}
Let $G_n=(V_n,E_n), n\in\N,$ be a sequence of finite graphs admitting a random weak limit $\rho$. Then,
\begin{equation}
\label{lim:g}\frac{\nu(G_n)}{|V_n|} \xrightarrow[n\to\infty]{} \gamma,
\end{equation}
where $\gamma\in [0,\frac{1}{2}]$ is defined by a recursion defined directly on the random weak limit $\rho$.
\end{theorem}

Since the work of Heilmann and Lieb \cite{heilmannlieb}, it is known
that the thermodynamic limit for monomer-dimer systems exists and
basic properties of this limit are derived for lattices. In
particular, \cite[Lemma 8.7]{heilmannlieb} shows the convergence of
the normalized matching number when the underlying graph is a lattice. More recently, Elek and Lippner \cite{eleklippner} extended this result by using the framework of local weak convergence for bounded degree graphs. Here we remove the bounded degree assumption. More importantly, the approach in \cite{eleklippner} is non-constructive. In contrast, we provide a characterization of $\gamma$ in terms of a local recursion defined directly on the random weak limit $\rho$. We postpone the discussion on how to actually compute $\gamma$ from $\rho$ to Subsection \ref{subsec:algo}. Our approach starts as in \cite{heilmannlieb} with the introduction of a natural family of probability distributions on the set of matchings parametrized by a single parameter $z>0$ called the Boltzmann-Gibbs distribution. The analysis in \cite{heilmannlieb} concentrates on the properties of the partition function, also known as the matching polynomial, from which a result like (\ref{lim:g}) can be deduced. Our analysis differs from this approach and concentrates on the analysis of the local marginals of the Boltzmann-Gibbs distribution, in a similar spirit as in the (non-rigorous) work of Zdeborov\'a and M\'ezard \cite{zdeborovamezard}. As in \cite{zdeborovamezard} and \cite{god81}, we start from an elementary formal recursion satisfied by the matching polynomials, and deduce an exact recursion for the local marginals of the Boltzmann-Gibbs distribution on any finite graph. A careful analysis of the contractivity properties of this recursion allows us to define the monomer-dimer model on infinite graphs (see Theorem \ref{co:cvweak}), and to define their "normalized matching number" (see Theorem \ref{co:monotony}). We should stress that the analysis of the marginal probabilities is essential for our second main result, namely the explicit computation of the matching number when the local weak limit is a Galton-Watson tree. Although simple adaptations of the argument in \cite{heilmannlieb} would yield a result like Theorem \ref{th:main}, the limit would be given in an implicit way which would not be sufficient to get our second main result.

As many other classical graph sequences, Erd\H{o}s-R\'enyi graphs and random graphs with a prescribed degree profile admit almost surely a particularly simple random weak limit, namely a \textit{unimodular Galton-Watson (UGW) tree} (see Example 1.1 in \cite{aldlyo}). This random rooted tree is parametrized by a probability distribution $\pi\in\cP(\N)$ with finite mean, called its \textit{degree distribution}. It is obtained by a Galton-Watson branching process where the root has offspring distribution $\pi$ and all other genitors have offspring distribution $\widehat{\pi}\in\cP(\N)$ defined by
\begin{equation}\label{eq:F}
\forall n\in \N, \widehat{\pi}_{n} = \frac{(n+1) \pi_{n+1}}{\sum_{k} k \pi_k}.
\end{equation}
Thanks to the markovian nature of the branching process, the recursion defining $\gamma$ simplifies into a \textit{recursive distributional equation}, which has been explicitly solved by the authors in a different context \cite{rankarxiv}. 
\begin{theorem}
\label{th:KS}
With the notation of Theorem \ref{th:main}, if the random weak limit $\rho$ is a UGW tree with degree distribution $\pi$, we have the explicit formula 
\begin{eqnarray*}
\gamma = \frac{1-\max_{t\in [0,1]} F(t)}{2},
\end{eqnarray*}
where 
\begin{eqnarray*}
F(t)=t\phi'(1-t)+\phi(1-t)+\phi\left(1-\frac{\phi'(1-t)}{\phi'(1)}\right)-1, 
\end{eqnarray*}
and $\phi(t)=\sum_{n}\pi_nt^n$ is the moment generating function of the degree distribution $\pi$. 
\end{theorem}

Differentiating the above expression, we see that any $t$ achieving the maximum must satisfy 
\begin{equation}
\label{eq:diff}
\phi'(1)t=\phi'\left(1-\frac{\phi'(1-t)}{\phi'(1)}\right).
\end{equation}
For Erd\H{o}s-R\'enyi random graphs with connectivity $c$, the degree
of the limiting UGW tree is Poisson with parameter $c$
(i.e. $\phi(t)=\exp(ct-c)$), so that (\ref{eq:diff}) becomes
$t=e^{-ce^{-ct}}$. We thus recover precisely Karp and Sipser's formula
(\ref{eq:KS}). Similarly, for random graphs with a prescribed degree
sequence, the log-concave assumption made by Bohmann and Frieze
guarantees that the above maximum is achieved at $t=0$ with $F(0)=0$, hence (\ref{eq:BF}) follows automatically.

A classical area of combinatorial optimization is formed by bipartite matching \cite{lovaszplummer}. 
We end this section, with a specialization of our results to bipartite graphs $G=(V=V^a\cup V^b, E)$.
The natural limit for a sequence of bipartite graphs is the following
hierarchal Galton-Watson tree  parameterized by two distributions
on $\N$ with finite first moment,
$\pi^a$ and $\pi^b$ and a parameter $\lambda\in [0,1]$.
We denote $\hpi^a$ and
$\hpi^b$ the corresponding distributions given by the transformation
(\ref{eq:F}).
We also denote $\phi^a$ and $\phi^b$ the generating functions of
$\pi^a$ and $\pi^b$.
The hierarchal Galton-Watson tree is then defined as follows: with probability $\lambda$, the
root has offspring distribution $\pi^a$, all odd generation genitors
have offspring distribution $\hpi^b$ and all even generation genitors
have offspring distribution $\hpi^a$; similarly with probability
$1-\lambda$, the root has offspring distribution $\pi^b$, all odd
generation genitors have offspring distribution $\hpi^a$ and all even
generation genitors have offspring distribution $\hpi^b$. In the first
(resp. second) case, we say that the root and all even generations are
of type $a$ (resp. $b$) and all the odd generations are of type $b$
(resp. $a$). To get a {\em unimodular hierarchal Galton-Watson (UHGW) tree} with degree distributions $\pi^a$ and $\pi^b$, we need to have:
$\lambda \phi^a{}'(1) =(1-\lambda)\phi^b{}'(1)$, so that 
\begin{eqnarray}
\label{eq:contrlamb}\lambda = \frac{\phi^b{}'(1)}{\phi^a{}'(1)+\phi^b{}'(1)}.
\end{eqnarray}
\begin{theorem}
\label{th:KSbi}
With the notation of Theorem \ref{th:main}, assume that the random weak limit $\rho$ is a UHGW tree with degree distributions
$\pi^a,\pi^b$.
If $\pi^a$ and $\pi^b$ have finite first moment, then
\begin{equation}
\label{eq:KSbi}
\gamma=\frac{\phi^b{}'(1)}{\phi^a{}'(1)+\phi^b{}'(1)}(1-\max_{t\in[0,1]}F^a(t)),
\end{equation}
where $F^a$ is defined by:
\begin{eqnarray*}
F^a(t) &=& \phi^a\left(1-\frac{\phi^b{}'(1-t)}{\phi^b{}'(1)}\right)-\frac{\phi^a{}'(1)}{\phi^b{}'(1)}\left(1-\phi^b(1-t)-t\phi^b{}'(1-t)\right).
\end{eqnarray*}
\end{theorem}
Note that if $\phi^a(x)=\phi^b(x)$, we find the same limit as in
Theorem \ref{th:KS}. Note that it is not obvious from formula
\eqref{eq:KSbi} that our expression for $\gamma$ is symmetric in $a$
and $b$ as it should. In the forthcoming Section \ref{sec:tlgwt},
Equation \eqref{eq:KSbi2} gives an alternative symmetric formula for
$\gamma$ which simplifies to (\ref{eq:KSbi}) thanks to (\ref{eq:contrlamb}).

Note also that our Theorem \ref{th:KSbi} computes the independence
number of random bipartite graphs.
Recall that a set of vertices in a graph $G$ is said to be
independent if no two of them are adjacent. The cardinality of any
largest independent set of points in $G$ is known as the independence
number of $G$ or the stability number of $G$ and is denoted by
$\alpha(G)$.
By K\H{o}nig's theorem, we know that for any bipartite graph
$G$ with vertex set $V$, $\alpha(G)+\nu(G)=|V|$.
The fact that a limit for $\frac{\alpha(G_n)}{{|V_n|}}$ exists,
has been proved recently in \cite{bayati-2010} for Erd\H{o}s-R\'enyi and
random regular graphs. The actual value for this limit is unknown except
for Erd\H{o}s-R\'enyi graphs with mean degree $c<e$. In this case, the
leaf-removal algorithm allows to compute explicitly the limit which agrees with
(\ref{eq:KSbi}) with $\phi^a(x)=\phi^b(x)=\exp(cx-x)$.

Motivated by some applications for Cuckoo Hashing \cite{fopa09},
\cite{kxor}, recent results have been obtained in the particular case
where $\pi^a(k)=1$ for some $k\geq 3$ and $\pi^b$ is a Poisson
distribution with parameter $\alpha k$. These degree distributions arise if one consider a sequence of bipartite graphs with $\lfloor \alpha m\rfloor$ nodes of type $a$ (called the items), $m$ nodes of type $b$ (called the locations) and each node of type $a$ is connected with $k$ nodes of type $b$ chosen uniformly at random (corresponding to the assigned locations the item can be stored in).
The result in this domain, obtained in \cite{frpa09} follows (see Section
\ref{sec:corfr}) from our Theorem \ref{th:KSbi}, namely:
\begin{corollary}\label{cor:fr}
Under the assumption of Theorem \ref{th:KSbi} and with $\pi^a(k)=1$ for some $k\geq 3$ and $\pi^b$ is a Poisson
distribution with parameter $\alpha k$.
Let $\xi$ be the unique solution of the equation:
\begin{eqnarray*}
k=\frac{\xi(1-e^{-\xi})}{1-e^{-\xi}-\xi e^{-\xi}},
\end{eqnarray*}
and $\alpha_c = \frac{\xi}{k(1-e^{-\xi})^{k-1}}$.
\begin{itemize}
\item for $\alpha\leq \alpha_c$, all (except $o_p(n)$) vertices of type
  $a$ are covered, i.e. $\frac{\nu(G_n)}{|V^a_n|} \xrightarrow[n\to\infty]{}1$.
\item for $\alpha>\alpha_c$, 
we have:
\begin{eqnarray}
\label{eq:min}\frac{\nu(G_n)}{|V^a_n|}
\xrightarrow[n\to\infty]{}1-\frac{1}{\alpha}\left(e^{-\xi^*}+\xi^*e^{-\xi^*}+\frac{\xi^*}{k}(1-e^{-\xi^*})-1\right),
\end{eqnarray}
where $\xi^*= k\alpha x^*$ and $x^*$ is the largest solution of
$x=\left(1-e^{-k\alpha x}\right)^{k-1}$.
\end{itemize}
\end{corollary}

In words, $\alpha_c$ is the load threshold: if $\alpha\leq \alpha_c$,
there is an assignment of the $\lfloor \alpha m\rfloor$ items to a
table with $m$ locations that respects the choices of all items,
whereas for $\alpha>\alpha_c$, such an assignement does not exist and
(\ref{eq:min}) gives the maximal number of items assigned without
collision. Note that results in \cite{fopa09}, \cite{kxor} are
slightly different in the sense that for the specific sequence of
random graphs described above (i.e. uniform hypergraphs), they show that for $\alpha<
\alpha_c$ all vertices of type $a$ are covered with high probability.
It is shown in \cite{lel12} how to get such results from Corollary
\ref{cor:fr} under the additional assumption that the sequence of
graphs are uniform hypergraphs.

\section{The Monomer-Dimer model} 
\label{sec:z}

We start with the case of a finite graph $G=(V,E)$.
Consider a natural family of probability
distributions on the set of matchings $\M(G)$,
parameterized by a single parameter $z>0$ called the
\textit{temperature} (note that the standard temperature $T$ in physics would correspond to $z=e^{-1/T}$ but this will not be important here): for any $M\in\M(G)$,
\begin{equation}
\label{eq:gibbs}
\mu^z_{G}(M)=\frac{z^{|V|-2|M|}}{P_G(z)},
\end{equation}
where $P_G$ is the \textit{matching polynomial},
$P_G(z)=\sum_{M\in{\M(G)}}z^{|V|-2|M|}$.
In statistical physics, this is called the \textit{monomer-dimer model} at temperature $z$ on $G$ (see \cite{heilmannlieb} for a complete treatment).
We let $\cM^z_G$ denote a random element of $\M(G)$ with
law $\mu^z_G$, and we call it a \textit{Boltzmann random matching at
  temperature $z$} on $G$.
Note that the lowest degree coefficient of $P_G$ is precisely the number of largest matchings on $G$. Therefore, $\mathcal M^z_G$ converges in law to a uniform largest matching as the temperature $z$ tends to zero.
We define the \textit{root-exposure probability (REP)} of the
rooted graph $(G,\circ)$ as
\begin{eqnarray}
\label{eq:mpe}\cR_z{[G,\circ]} = \mu^z_{G}\left(\circ\textrm{ is
    exposed}\right).
\end{eqnarray}

\subsection{Local recursions to compute $\gamma$}\label{subsec:algo}


Before starting with the proof, we explain (whithout proofs) how to
compute $\gamma$ in (\ref{lim:g}).
For a finite graph, our computation of $\gamma$ follows exactly the approach
of Godsil \cite{god81}.
We recall Godsil's notion of
the \textit{path-tree} associated with a rooted graph $G$:
if $G$ is any rooted graph with root $\circ$, we define
its path-tree $T_G$ as the rooted tree whose vertex-set consists of all finite simple paths starting at the
  root $\circ$; whose edges are the pairs $\{P,P'\}$ of the
  form $P=v_1\ldots v_n$, $P'=v_1\ldots v_nv_{n+1} (n\geq 1)$; whose root is the single-vertex path $\circ$.
By a \textit{finite simple path}, we mean here a finite sequence of
distinct vertices $v_1\ldots v_n$ ($n\geq 1$) such that  $v_iv_{i+1}\in E$
for all $1\leq i < n$.

It is well-known since Godsil's result \cite{god81} that path-trees
capture considerable information about matchings in general graph and
are easier to work with than the graph itself.
For a rooted graph $[G,\circ]$, let $T_{[G,\circ]}$ be the associated
path-tree and consider the corresponding system of equations (where $u\succ
v$ if $u$ is a child of $v$):
\begin{equation}
\label{eq:fpzeroi}
\forall v\in T_{[G,\circ]}, \quad x_v=\frac {1}
{1+\sum_{u\succ v}{\left({\sum_{w\succ u}x_w}\right)^{-1}}}.
\end{equation}
For any finite rooted graph $[G,\circ]$, (\ref{eq:fpzeroi}) has a
unique solution in $[0,1]^{T_{[G,\circ]}}$ and we denote the value taken at the
root by $x_\circ(G)$. Then $x_\circ(G)$ is exactly the probability for
the root $\circ$ of being exposed in a uniform maximal matching. In
particular, we have
\begin{eqnarray*}
\nu(G) = \sum_{v\in V}\frac{1-x_v(G)}{2}.
\end{eqnarray*}
This argument follows from \cite{god81} and will be a special case of our
analysis.

For an infinite graph with bounded degree, it turns out that it is not
always possible to make sense of the local recursions
(\ref{eq:fpzeroi}). However, our analysis will show that for any $z>0$, the
infinite set of equations: 
\begin{equation*}
\forall v\in T_{[G,\circ]}, \quad x_v(z)=\frac {1}
{1+\sum_{u\succ v}{\left(z^2+{\sum_{w\succ u}x_w(z)}\right)^{-1}}},
\end{equation*}
has a unique solution in $[0,1]^{T_{[G,\circ]}}$ and the value taken
by the root is exactly $\cR_z[G,\circ]$ (which is the probability for the root
$\circ$ of being exposed in a Boltzmann random matching at temperature
$z$ when the graph $G$ is finite). Then our Theorem \ref{co:cvMN} will
imply that for any sequence of
finite graphs $(G_n=(V_n,E_n))_{n\in\N}$ satisfying $|E_n|=O(|V_n|)$ and having $\rho$ as a random weak limit,
\begin{eqnarray*}
\frac{\nu(G_n)}{{|V_n|}}\xrightarrow[n\to\infty]{}\frac{1-\EE_\rho\left[\lim_{z\to
      0}\cR_z\right]}{2},
\end{eqnarray*}
and $\lim_{z\to 0}\cR_z[G,\circ]$ is actually the largest solution to (\ref{eq:fpzeroi}).
From a practical point of view, it is possible to compute an
approximation of $\cR_z[G,\circ]$ by looking at a sufficient large
ball centered at the root $\circ$. 
Moreover our analysis will show that the quantity
$\EE_\rho\left[\cR_z\right]$ is a good approximation of $\EE_\rho\left[\lim_{z\to
      0}\cR_z\right]$ as soon as $|E_n|=O(|V_n|)$ (see Lemma \ref{lm:unifctrl}).

\subsection{Extension of the model on infinite graphs with bounded degree}
Let $G-\circ$ be the graph obtained from $G$ by removing its root $\circ$.
Since the matchings of $G$ that expose $\circ$ are exactly the matchings of
$G-\circ$, we have the identity
\begin{equation}
\cR_z{[G,\circ]}= \frac{z P_{G-\circ}(z)}{P_G(z)},
\end{equation}
which already shows that the REP is an analytic function of the
temperature. The remarkable fact that its domain of analyticity contains the
right complex half-plane $$\HH_+=\{z\in\C;\Re(z)>0\}$$
is a consequence of the powerful Heilmann-Lieb theorem \cite[Theorem 4.2]{heilmannlieb}
(see \cite{sokalwagner} for generalizations).
The key to the study of the REP is the following elementary but fundamental local recursion :
\begin{equation}
\label{eq:recg}
\cR_z{[G,\circ]}=z^2\left(z^2+\sum_{v\sim \circ}\cR_z{[G-\circ,v]}\right)^{-1}.
\end{equation}
Clearly, this recursion determines uniquely the functional $\cR_z$ on the
class of finite rooted graphs, and may thus be viewed as an inductive
definition of the REP.  Remarkably enough, this alternative characterization
allows for a continuous extension to infinite graphs with bounded degree, even though the above recursion never ends. We let $\cH$ denote the space of analytic functions on $\HH_+$,  equipped with its usual topology of uniform
convergence on compact sets. Our fundamental lemma is as follows :
\begin{theorem}[The fundamental local lemma]\mbox{}
\label{th:mainb}
\begin{enumerate}
\item For every fixed $z\in\HH_+$,  the local recursion (\ref{eq:recg})
  determines a unique $\cR_z\colon \cGr \to z\HH_+$.
\item For every fixed $[G,\circ]\in\cGr$, $z\mapsto \cR_{z}[G,\circ]$ is
  analytic.
\item The resulting mapping $[G,\circ] \in\cGr\longmapsto
  \cR_{(\cdot)}[G,\circ]\in\cH$ is continuous.
\end{enumerate}
\end{theorem}
This local lemma has strong implications for the monomer-dimer model, which
we now list. The first one is the existence of an infinite volume limit for
the Gibbs-Boltzmann distribution.

\begin{theorem}[Monomer-dimer model on infinite graphs]
\label{co:cvweak}
Consider a graph $G\in\cGr$ and a temperature $z>0$. For any finite matching $M$ of $G$, the cylinder-event marginals defined by
\begin{equation*}
\mu_G^z(M\subseteq \cM) = z^{-2|M|}\prod_{k=1}^{2|M|}\cR_z{[G - \{v_1,\ldots,v_{k-1}\},v_k]},
\end{equation*}
are consistent and independent of the ordering $v_1,\ldots v_{2|M|}$ of the
vertices spanned by $M$. They thus determine a unique probability
distribution $\mu_G^z$ over the matchings of $G$. It coincides with the former definition in the case where $G$ is
finite, and extends it continuously in the following sense : for any $\circ\in V$ and any sequence $([G_n,\circ_n])_{n\in\N}\in\cGr^\N$ converging to $[G,\circ]$,
$$[G_n,\circ_n,\cM_n]\xrightarrow[n\to\infty]{d}[G,\circ,\cM],$$
in the local weak sense for random networks, where $\cM_n$ has law $\mu^z_{G_n}$ and $\cM$ has law $\mu_G^z$.
\end{theorem}

%
Although it is not our concern here, we obtain as a by-product the strong convergence of the logarithm of the matching polynomial, also called \textit{free energy} in the monomer-dimer model :
\begin{corollary}
\label{co:fe}
Let $(G_n)_{n\in\N}$ be a sequence of finite graphs with bounded degree admitting a random weak limit $[G,\circ]$. The following convergence holds in the analytic sense on $\HH_+$,
\begin{equation*}
\frac{1}{|V_n|}\log \frac{P_{G_n}(z)}{P_{G_n}(1)}\xrightarrow[n\to\infty]{}\int_1^{z}\frac{\EE_{\rho}[\mathcal{R}_s[G,\circ]]}{s}ds,
\end{equation*}
where $\EE_\rho[\mathcal{R}_s[G,\circ]]$ denotes  the expectation under the measure $\rho$ of the variable $\mathcal{R}_s[G,\circ]$.
\end{corollary}
A similar result was established in \cite{heilmannlieb} for the lattice case, and in \cite{bayatinair} under a restrictive large girth assumption. 

\subsection{Proof of Theorem \ref{th:mainb} : the fundamental lemma}
\label{sec:finite}

The local recursion (\ref{eq:recg}) involves mappings of the form :
$$\phi_{z,d}\colon \left(x_1,\ldots,x_d\right) \mapsto z^2\left(z^2+\sum_{i=1}^d
  x_i\right)^{-1},$$
where $d\in\N$. In the following lemma, we gather a few elementary
properties of this transformation, which are immediate to check but will be of constant use throughout
the paper.

\begin{lemma}[Elementary properties]
\label{lm:propphi}
For any $d\in\N$ and $z\in\HH_+$,
\begin{enumerate}
\item $\phi_{z,d}$ maps analytically $z\HH_+\times\ldots\times z\HH_+$ into $z\HH_+$
\item $|\phi_{z,d}|$ is uniformly bounded by $|z|/\Re(z)$ on
  $z\HH_+\times\ldots\times z\HH_+$.
\end{enumerate}
\end{lemma}
From part $1$, it follows that the
REP of a finite rooted graph belongs to $\cH$, when viewed as a function of the temperature $z$. Part $2$ and Montel's theorem guarantee that the family of all those REPs
is relatively compact in $\cH$.  Note that relative compactness also plays a central role in \cite{heilmannlieb}. Combined with the following
uniqueness property at high temperature, it will quickly lead to the proof of
Theorem \ref{th:mainb}.

The local recursion (\ref{eq:recg}) also involves graph
transformations of the form $[G,\circ]\mapsto[G-\circ,v]$, where
$v\sim\circ$. Starting from a given $[G,\circ]\in\cGr$,  we let
$\rm{Succ}^*[G,\circ]\subseteq\cGr$ denote the (denumerable) set of all
rooted graphs that can be obtained by successively applying finitely many
such transformations.  
\begin{lemma}[Uniqueness at high temperature]
\label{lm:uniqueness}
Let $[G,\circ]\in\cGr$ and $z\in \HH_+$ such that  $\Re(z)> \sqrt{ d}$. If
$$\cR^1_z,\cR_z^2\colon \rm{Succ}^*[G,\circ]\to z\HH_+$$ both satisfy the local recursion
(\ref{eq:recg}) then $\cR^1_z=\cR^2_z$.
\end{lemma}

\bp{}Set $\alpha=2{|z|}/{\Re(z)}$ and $\beta=\Re(z)^{-2}$. From
(\ref{eq:recg}) and part 2 of Lemma \ref{lm:propphi} it is clear that the absolute difference
$\Delta=|\cR^1_z-\cR^2_z|$ must satisfy
$$\Delta[G,\circ]\leq \alpha\qquad\textrm{ and }\qquad\Delta[G,\circ]\leq \beta\sum_{v\sim\circ}\Delta[G-\circ,v].$$
In turn, each $\Delta[G-\circ,v]$ appearing in the second upper-bound may
be further expanded into
$\beta\sum_{w\sim v,w\neq\circ}\Delta[G-\circ-v,w]$. Iterating this procedure $k$ times, one obtains 
$\Delta[G,\circ]\leq \beta^{k}d^k\alpha$.
Taking the infimum over all $k$
yields $\Delta[G,\circ]=0$, since the assumption $\Re (z ) > \sqrt d$ means precisely $\beta d<1$.
\ep

\bp{ of Theorem \ref{th:mainb}} For clarity, we divide the proof in three parts : we first define a specific solution  which satisfies (\ref{eq:recg}). We will then prove its unicity and check its continuity. This will prove parts 1-3 of Theorem \ref{th:mainb}.

\textbf{Analytic existence. }
Fix $[G,\circ]\in\cGr$, and consider an
arbitrary collection of $\HH_+\to z\HH_+$ analytic functions
$z\mapsto \cR^0_z[H,i]$, indexed by the elements $[H,i]\in\rm{Succ}^*[G,\circ]$. For every $n\geq 1$, define recursively
\begin{equation}
\label{eq:def}
\cR^{n}_z[H,i]=z^2\left(z^2+\sum_{j\sim i}\cR^{n-1}_{z}[H-i,j]\right)^{-1},
\end{equation}
for all $z\in\HH_+$ and $[H,i]\in\rm{Succ}^*[G,\circ]$. By Lemma
\ref{lm:propphi}, each sequence  $\left(z\mapsto \cR^n_{z}[H,i]\right)_{n\in\N}$
is relatively compact in $\cH$. Consequently, their joint collection
as $[H,i]$ varies in the denumerable set $\rm{Succ^*}[G,\circ]$ is
relatively compact in the product space
$\cH^{\rm{Succ^*}[G,\circ]}$. Passing to the limit in (\ref{eq:def}), we
see that any pre-limit $\cR_z\colon\rm{Succ^*}[G,\circ]\to z\HH_+$ must automatically satisfy (\ref{eq:recg})
for each $z\in\HH_+$. By Lemma \ref{lm:uniqueness}, this determines uniquely the
value of $\cR_z[G,\circ]$ for $z$ with sufficiently large real part, and hence
everywhere in $\HH_+$ by analyticity. To sum up, we have just proved the following : for every $[G,\circ]\in\cGr$,
the limit
\begin{equation}
\label{eq:recsol}
\cR_z{[G,\circ]}:=\lim_{n\to\infty}\cR^n_z{[G,\circ]}
\end{equation}
exists in $\cH$, satisfies the recursion (\ref{eq:recg}), and does not
depend upon the choice of the initial condition $\cR^0_z\colon
\rm{Succ^*}[G,\circ]\to z\HH_+$ (provided that the latter is analytic
in $z\in\HH_+$).

\textbf{Pointwise uniqueness. }Let us now show that any
$\cS\colon\rm{Succ^*}[G,\circ]\to z\HH_+$ satisfying the recursion (\ref{eq:recg}) at a fixed value $z=z_0\in\HH_+$ must coincide with the
$z=z_0$ specialization of the analytic solution constructed above. For
each $[H,i]\in\rm{Succ^*}[G,\circ]$, the constant initial function $\cR^{0}_{z}[H,i]:=\cS[H,i]$ is trivially analytic from $\HH_+$ to $z\HH_+$,
so the iteration (\ref{eq:def}) must converge to the analytic solution
$\cR_z$. Since $\cR^n_{z_0}=\cS$ for all $n\in\N$, we obtain $\cR_{z_0}=
\cS$, as desired.

\textbf{Continuity. }Finally, assume that $\left([G_n,\circ]\right)_{n\geq
  1}\in\cGr^\N$ converges locally to $[G,\circ]$, and let us show that
\begin{equation}
\label{eq:continuity}
\cR_z{[G_n,\circ]}\xrightarrow[n\to\infty]{\cH}\cR_z[G,\circ].
\end{equation}
It is routine that,
up to rooted isomorphisms, $G,G_1,G_2,\ldots$ may be represented on a common vertex set, in such a
way that for each fixed $k\in\N$, $[G_n,\circ]_k=[G,\circ]_k$ for all $n\geq n_k$. By construction, any
simple path $v_1\ldots v_{k}$
starting from the root in $G$ is now also a simple path starting from the root
in each $G_n,n\geq n_k$, so the $\cH-$valued sequence
$\left(z\mapsto \cR_{z}[G_n-\{v_1,\ldots,v_{k-1}\},v_k]\right)_{n\geq n_k}$ is
well defined, and relatively compact (Lemma \ref{lm:propphi}). Again, the denumerable
collection of all sequences obtained by letting the simple path
$v_1 \ldots v_{k}$ vary in $[G,\circ]$ is relatively compact for the product topology,
and any pre-limit must by construction satisfy (\ref{eq:recg}). By
pointwise uniqueness, the convergence (\ref{eq:continuity}) must hold. \ep

\subsection{Proof of Theorem \ref{co:cvweak} : convergence of the Boltzmann distribution}
\label{sec:cvweak}
Consider an infinite $[G,\circ]\in\cGr$, and let
$\left([G_n,\circ]\right)_{n\geq 1}$ be a sequence of finite rooted
connected graphs converging locally to $[G,\circ]$. As above, represent $G,G_1,G_2,\ldots$ on a common vertex set, in such a
way that for each $k\in\N$, $[G_n,\circ]_k=[G,\circ]_k$ for all $n\geq
n_k$.
Now fix an arbitrary finite matching $M$ in $G$, and denote by
$v_1,\ldots,v_{2|M|}$ the vertices spanned by $M$, in any order. By
construction, $M$ is also a matching of $G_n$ for large enough $n$. But the matchings of $G_n$ that contain $M$ are in $1-1$ correspondence with the matchings of $G_n -
\{v_1,\ldots,v_{2|M|}\}$, and hence
$$\mu^z_{G_n}\left(M\subseteq\cM\right)=\frac{P_{G_n-\{v_1,\ldots,v_{2|M|}\}}(z)}{P_{G_n}(z)}=z^{-2M}\prod_{k=1}^{2M}\cR_z[G_n-\{v_1,\ldots,v_{k-1}\},v_k].$$
But  $[G_n-\{v_1,\ldots,v_{k-1}\},v_k]$ converges locally to
$[G-\{v_1,\ldots,v_{k-1}\},v_k]$, so by continuity of $\cR_z$,
$$\mu^z_{G_n}\left(M\subseteq\cM\right)\xrightarrow[n\to\infty]{}z^{-2M}\prod_{k=1}^{2M}\cR_z[G-\{v_1,\ldots,v_{k-1}\},v_k].$$

\bp{ of Corollary \ref{co:fe}}
Analytic convergence of the free energy follows from Theorem \ref{co:cvweak} and Lebesgue dominated convergence Theorem, since for any finite graph $G=(V,E)$ we have
$$(\log P_G)'(z) = \frac{P_G'(z)}{P_G(z)} = \frac{1}{|V|} \sum_{\circ \in V}\frac{\cR_z[G,\circ]}{z} = \frac{\rho[\cR_z[G,\circ]]}{z}.$$
The uniform domination $\left|\frac{\rho[\cR_z[G,\circ]]}{z}\right|\leq \frac{1}{\Re(z)}$ is provided by Lemma \ref{lm:propphi}.
\ep

\section{The zero-temperature limit}\label{sec:zero}

Motivated by the asymptotic study of maximum matchings, we now let the
temperature $z\to 0$.

\subsection{The case of bounded degree}

We first use the results from the previous section to prove a version of Theorem \ref{th:main} for graphs with bounded degree. 
\begin{theorem}[The zero temperature limit in graphs with bounded degree]
\label{co:monotony}
For any $[G,\circ]\in\cGr$, the zero temperature limit
$$\cR_*[G,\circ]=\lim_{z\to 0}\downarrow \cR_z[G,\circ]$$
exists. Moreover, $\cR_*\colon \cGr \to [0,1]$
is the largest solution to the recursion
\begin{equation}
\label{eq:recgi}
\cR_*[G,\circ] = \left(1+\sum_{v\sim\circ}\left(\sum_{w\sim
     v}\cR_*[G-\circ-v,w]\right)^{-1}\right)^{-1},
\end{equation}
with the conventions $0^{-1}=\infty$, $\infty^{-1}=0$. When $G$ is
finite,  $\cR_*[G,\circ]$ is the probability that $\circ$ is exposed
in a uniform maximum matching.
\end{theorem}
\bp{}
Fix $[G,\circ]\in\cGr$. First, we claim that $z\mapsto \cR_z[G,\circ]$ is
non-decreasing on $\R_+$. Indeed, this is obvious if $G$ is reduced to
$\circ$, since in that case the REP
is simply $1$. It then inductively extends to any finite graph $[G,\circ]$,
because iterating twice (\ref{eq:recg}) gives
\begin{equation}
\label{eq:recgtwice}
\cR_z[G,\circ] = \left(1+\sum_{v\sim\circ}\left(z^2+\sum_{w\sim
     v}\cR_z[G-\circ-v,w]\right)^{-1}\right)^{-1}.
\end{equation}
For the infinite case, $[G,\circ]$ is the local limit of the sequence
of finite truncations $\left([G,\circ]_n\right)_{n\in\N}$, where, for $n \geq 1$, $[G,\circ]_n$ denotes the finite rooted subgraph induced by the vertices lying at graph-distance at
most $n$ from $\circ$. So by
continuity of the REP, $\cR_z[G,\circ]=\lim_{n\to\infty}\cR_z[G,\circ]_n$
must be non-decreasing in $z$ as well.
This guarantees the existence of the $[0,1]-$valued limit
$$\cR_*[G,\circ]=\lim_{z\to 0}\downarrow\cR_z[G,\circ].$$
Moreover, taking the $z\to 0$ limit in (\ref{eq:recgtwice}) guarantees
the recursive formula (\ref{eq:recgi}).

Finally, consider
$\cS_*\colon\rm{Succ}^*[G,\circ]\to[0,1]$ satisfying the recursion
(\ref{eq:recgi}). Let us show by induction over $n\in\N$ that for every
$[H,i]\in\rm{Succ}^*[G,\circ]$ and $z>0$,
\begin{equation}
\label{eq:compare}
\cS_*[H,i]\leq \cR_z[H,i]_{2n}.
\end{equation}
The statement is trivial when $n=0$ ($\cR_z[H,i]_0=1$), and is preserved from $n$ to $n+1$
because
\begin{eqnarray*}
\cR_z[H,i]_{2n+2} & = & \left(1+\sum_{j\sim i}\left(z^{2}+\sum_{k\sim
      j}\cR_z[H-i-j,k]_{2n}\right)^{-1}\right)^{-1} \\
& \geq & \left(1+\sum_{j\sim
    i}\left(\sum_{k\sim j}\cS_*[H-i-j,k]\right)^{-1}\right)^{-1} =\, \cS_*[H-i,j].
\end{eqnarray*}
Letting $n\to\infty$ and then $z\to 0$ in (\ref{eq:compare}) yields $\cS_*\leq \cR_*$, which completes the proof
\ep

This naturally raises the following question : may the
\textit{zero temperature limit} be interchanged with the \textit{infinite
  volume limit}, as suggested by the diagram below ?

\begin{figure}[h!]
$$\xymatrix @!0 @R=2.5cm @C=4cm{
   \cR_z[G_n,\circ_n]\ar[d]_{z\to 0 } \ar[r]^{n\to\infty}
   & \cR_z[G,\circ]\ar[d]^{z\to 0} \\
   \cR_*[G_n,\circ_n] \ar@{.>}[r]_{n\to\infty} & \cR_*[G,\circ]
 }$$
\label{fig:proof}
\end{figure}

Unfortunately, the recursion (\ref{eq:recgi}) may admit several distinct
solutions, and this translates as follows : in the limit of zero temperature, \textit{correlation decay breaks for the monomer-dimer model}, in the precise sense that
the functional $\cR_*\colon\cGr\to[0,1]$ is no longer continuous with respect to local convergence. For example, one can easily construct an infinite rooted tree $[T,\circ]$ with
bounded degree such that
\begin{eqnarray*}
\lim_{n\to\infty}\downarrow \cR_*[T,\circ]_{2n}  & \neq &
\lim_{n\to\infty}\uparrow \cR_*[T,\circ]_{2n+1}.
\end{eqnarray*}
Indeed, consider the case of $T$ being the graph on $\N$ rooted at $0 = \circ$, where two integers share an edge if they differ by $1$. Then, a straightforward computation gives $\cR_*[T,\circ]_{2n} =1/2$ while $\cR_*[T,\circ]_{2n+1}  = 0$. 
Despite this lack of correlation decay,
the interchange of limits turns out to be valid ``on
average'', i.e. when looking at a uniformly chosen vertex
$\circ$. 

\begin{theorem}[The limiting matching number of bounded-degree graph sequences]
\label{co:cvMN}
Let $\rho$ be a probability distribution over $\cGr$. For any sequence of
finite graphs $(G_n=(V_n,E_n))_{n\in\N}$ satisfying $|E_n|=O(|V_n|)$ and having $\rho$ as a random weak limit,
\begin{eqnarray*}
\frac{\nu(G_n)}{{|V_n|}}\xrightarrow[n\to\infty]{}\frac{1-\EE_\rho\left[\cR_*\right]}{2}.
\end{eqnarray*}
\end{theorem}

In order to get our Theorem \ref{th:main}, we need to remove the bounded degree assumption. This is done below.
In the case where the limit $\rho$ is a (two-level) Galton-Watson tree, the recursion
(\ref{eq:recgi}) simplifies into a \textit{recursive distributional
  equation} (RDE). The computations for these cases are done in Section \ref{sec:rdet}.

\bp{ of Theorem \ref{co:cvMN}}

Let $G=(V,E)$ be a finite graph and $M$ be any maximal matching of $G$. Then 
$$
\sum_{ v \in V} \ind (v \hbox{ is exposed in $M$}) =  |V| - 2 \sum_{ e \in E} \ind (e \in M). 
$$
In particular, if $\rho=\cU(G)$, we have the
elementary identity 
\begin{equation}\label{eq:R*nu}
\EE_\rho\left[\cR_*\right] = 1-\frac{2\nu(G)}{|V|}.
\end{equation} The proof of Theorem \ref{co:cvMN} will easily follow from the following
uniform control:
\begin{lemma}[Uniform continuity around the zero-temperature point]
\label{lm:unifctrl}
Let $G=(V,E)$ be a finite graph. For any $0<z<1$,
\begin{eqnarray}
\label{eq:lemma}
\EE_\rho\left[\cR_z\right] + \frac{|E|}{|V|}\frac{\log 2}{\log z} \leq
\EE_\rho\left[\cR_*\right] \leq \EE_\rho\left[\cR_z\right].
\end{eqnarray}
\end{lemma}
Indeed, let $\rho$ be a probability distribution on $\cGr$, and let
$(G_n=(V_n,E_n))_{n\in\N}$ be a sequence of finite graphs with
$|V_n|=O(|E_n|)$, whose random weak limit is $\rho$. For each $n\in\N$, set
$\rho_n=\cU(G_n)$. With these notations, proving Theorem \ref{co:cvMN}
amounts to establish :
\begin{equation}
\label{eq:cvrepi}
\EE_{\rho_n}\left[\cR_*\right]\xrightarrow[n\to\infty]{}\EE_\rho\left[\cR_*\right].
\end{equation}
However, since $\rho_n\Longrightarrow \rho$, and since each $\cR_z,z>0$ is continuous and
bounded, we have for every $z>0$,
$$\EE_{\rho_n}[\cR_z] \xrightarrow[n\to\infty]{} \EE_\rho[\cR_z].$$
Thus, setting  $C=\sup_{n\in\N}\frac{|E_n|}{|V_n|}$ and letting $n\to\infty$ in (\ref{eq:lemma}), we see that for any $z<1$,
\begin{eqnarray*}
\EE_\rho\left[\cR_z\right] + C\frac{\log 2}{\log z} \leq \liminf_{n\to\infty}\EE_{\rho_n}\left[\cR_*\right]  \leq   \limsup_{n\to\infty}\EE_{\rho_n}\left[\cR_*\right]  \leq \EE_{\rho}\left[\cR_z\right].
\end{eqnarray*}
Letting finally $z\to 0$, we obtain exactly (\ref{eq:cvrepi}), and it
only remains to show Lemma \ref{lm:unifctrl}.

\bp{ of Lemma \ref{lm:unifctrl}} Fix $0<z<1$. Since $z\mapsto
\EE_\rho\left[\cR_z\right]$ is non-decreasing, we have
$$ \EE_\rho\left[\cR_*\right] \leq \EE_\rho\left[\cR_z\right] \leq \frac{-1}{\log z}\int_{z}^1s^{-1} \EE_\rho\left[\cR_s\right]ds.$$
Use $\EE_\rho\left[\cR_s\right]= \frac{s P'_G(s)}{|V| P_G(s)}$ to rewrite this as
$$\EE_\rho\left[\cR_*\right] \leq \EE_\rho\left[\cR_z\right] \leq  \frac{1}{|V|\log z}\log
\frac{P_G(z)}{P_G(1)}.$$
Now, $P_G(1)$ is the total number of matchings and is thus clearly at
most $2^{|E|}$, while $P_G(z)$ is at least $z^{|V|-2\nu(G)}$. Using \eqref{eq:R*nu}, these two bounds yield to 
$$
\EE_\rho\left[\cR_*\right] \leq \EE_\rho\left[\cR_z\right] \leq  \frac{1}{|V|\log z} \left( |V|\EE_\rho\left[\cR_*\right] \log z  - |E| \log 2 \right).
$$
This gives (\ref{eq:lemma}). \ep

\subsection{The case of unbounded degree}
\label{sec:deg}
In this section, we establish Theorem \ref{th:main} in full generality, removing the restriction of bounded degree from Theorem \ref{co:cvMN}. To this end, we introduce the $d-$\textit{truncation} $G^d$ ($d\in \N$) of a graph $G=(V,E)$, obtained from $G$ by \textit{isolating} all vertices with degree more than $d$, i.e. removing any edge incident to them. This transformation is clearly continuous with respect to local convergence. Moreover, its effect on the matching number can be easily controlled : 
\begin{equation}
\label{eq:trunc}
\nu(G^d)\leq \nu(G)\leq \nu(G^d)+\#\{v\in V;\deg_G(v)>d\}.
\end{equation}
Now, consider a sequence of finite graphs $(G_n)_{n\in\N}$ admitting a random weak limit $(G,\circ)$. 
First, fixing $d\in\N$, we may apply Theorem \ref{co:cvMN} to the sequence $(G^d_n)_{n\in\N}$ to obtain :
\begin{equation*}
\frac{\nu(G^d_n)}{|V_n|}\xrightarrow[n\to\infty]{}\frac{1-\EE_{\rho_d}\left[\cR_*\right]}2,
\end{equation*}
where $\rho_d$ is the $d$-truncation of $\rho$.
Second, we may rewrite (\ref{eq:trunc}) as
\begin{equation*}
\left|\frac{\nu(G^d_n)}{|V_n|}-\frac{\nu(G_n)}{|V_n|}\right|\leq \frac{\#\{v\in V_n;\deg_{G_n}(v)>d\}}{|V_n|}.
\end{equation*}
Letting $n\to\infty$, we obtain
\begin{equation*}
\limsup_{n\to\infty}\left|\frac{1-\rho_d\left[\cR_*\right]}2-\frac{\nu(G_n)}{|V_n|}\right|\leq \rho\left(\deg(\circ)>d\right),
\end{equation*}
This last line is, by an elementary application of Cauchy criterion, enough to guarantee the convergence promised by Theorem \ref{th:main}, i.e.
\begin{equation}
\label{eq:toshow}
\frac{\nu(G_n)}{|V_n|}\xrightarrow[n\to\infty]{}\gamma,\qquad\textrm{ where }\qquad\gamma:=\lim_{d\to\infty}\frac{1-\EE_{\rho_d}\left[\cR_*\right]}2.
\end{equation}

Note that because of the possible absence of correlation decay, the largest solution $\cR_*[G,\circ]$ is not a continuous function of $(G,\circ)\in\cG$. In particular, we do not know whether it is always the case that 
\begin{equation}
\label{eq:question}
\gamma=\frac{1-\EE_\rho[\cR_*]}2,
\end{equation} as established in Theorem \ref{co:cvMN} for graphs with bounded degree. However, (\ref{eq:question}) holds in the particular cases where we have an explicit formula for $\EE_\rho[\cR_*]$ which depends continuously upon the degree distribution as will be the case in Section \ref{sec:rdet}. 

\section{Computations on (hierarchal) Galton-Watson trees}\label{sec:rdet}
\subsection{The case of Galton-Watson trees}
\label{sec:gwt}
We now investigate the special case where the limiting random graph is
a UGW tree $T$. Specifically, we fix a distribution
$\pi\in\cP(\N)$ with finite support (we will relax this
assumption in the sequel) and we consider
 a UGW tree $T$ with degree
 distribution $\pi$ as defined in Section \ref{sec:res}.
The random matchings $\cM^z_T, z\geq 0$ are perfectly
well-defined, and all the previously established results for graphs with bounded degree hold almost
surely. However, the self-similar
recursive structure of $T$ gives to the fixed-point characterizations
(\ref{eq:recg}) and (\ref{eq:recgi}) a very special form that is worth
making explicit.

Before we
start, let us insist on the fact that $\cR_z[T] $ ($z >0$) is random : it
is the quenched probability that the root is exposed at temperature $z$,
given the random tree $T$. In light of Theorem \ref{th:main}, it becomes important to ask for its distribution.
Let $\cP\left([0,1]\right)$ denote the space of Borel probability measures on $[0,1]$. Given $z>0$, $\nu\in\cP\left(\N\right)$ and $\mu\in\cP\left([0,1]\right)$, we denote by $\Theta_{\nu,z}(\mu)$ the law of the $[0,1]-$valued r.v.
\begin{equation*}
Y{=}\frac {z^2}
{z^2+\sum_{i=1}^{{\cN}}{X}_i},
\end{equation*}
where $\cN\sim \nu$ and ${X}_1,{X}_2,\ldots \sim \mu$, all of them being
independent. This defines an operator $\Theta_{\nu,z}$ on
$\cP\left([0,1]\right)$. The corresponding fixed point equation
$\mu=\Theta_{\nu,z}(\mu)$ belongs to the general class of \textit{recursive distributional equations}, or RDE.
Equivalently, it can be rewritten as
\begin{equation*}
\label{rde}
X\stackrel{d}{=}\frac {z^2}
{z^2+\sum_{i=1}^{{\cN}}{X}_i},
\end{equation*}
where ${X}_1,{X}_2,\ldots$ are i.i.d. copies of the unknown random variable ${X}$. Note that the same RDE appears in the analysis of the spectrum and rank of adjacency matrices of random graphs \cite{res10}, \cite{rankarxiv}.
With this notations in hands, the infinite system of equations
(\ref{eq:recg}) defining $\cR_z[T]$ clearly leads to the following
distributional characterization:
\begin{lemma}
For any $z>0$, $\cR_z[T]$ has distribution $\Theta_{\pi,z}(\mu_z)$,
where $\mu_z$ is solution to the RDE $\mu_z=\Theta_{\widehat{\pi},z}(\mu_z)$.
\end{lemma}

The same program can be carried out in the zero temperature limit. Specifically, given $\nu,{\nu'}\in\cP(\N)$ and $\mu\in\cP\left([0,1]\right)$, we
define $\Theta_{\nu,{\nu'}}(\mu)$ as the law of the $[0,1]-$valued r.v.
\begin{equation}
\label{eq:operator}
Y=\frac{1}{1+\sum_{i=1}^{\cN}\left(\sum_{j=1}^{{\cN_i}'}X_{ij}\right)^{-1}},
\end{equation}
where ${\cN}\sim\nu$, ${\cN_i}'\sim{\nu'}$, and ${X}_{ij}\sim\mu$, all of
them being independent. This defines an operator $\Theta_{\nu,{\nu'}}$ on
$\cP\left([0,1]\right)$ whose fixed points will play a crucial role in our
study. 
Then, Theorem \ref{co:monotony} implies:
\begin{lemma}
\label{th:mtzero}
The random variable $\cR_*[T]$ has law
$\Theta_{\pi,\widehat{\pi}}(\mu_*)$, where $\mu_*$ is the largest
solution to the RDE $\mu_*=\Theta_{\widehat{\pi},\widehat{\pi}}(\mu_*)$.
\end{lemma}
Recall that the mean of $\Theta_{\pi,\widehat{\pi}}(\mu_*)$
 gives precisely the asymptotic size of a maximum matching for
any sequence of finite random graphs whose random weak limit is
$T$ (Theorem \ref{co:cvMN}). We will solve this RDE in the next section in the more general set-up of UHGW trees. Combined with Theorem \ref{th:main} and a simple continuity argument to remove the bounded degree assumption, this will prove Theorem \ref{th:KS}.

\subsection{The case of hierarchal Galton-Watson trees}
\label{sec:tlgwt}

As in previous section, we first assume that both $\pi^a$ and $\pi^b$
have a finite support.
We can define a RDE but with some care
about the types $a$ and $b$. The corresponding results read as
follows:
\begin{lemma}
\label{lem:RDEbi}
For any $z>0$, conditionally on the root being of type $b$
(resp. $a$), $\cR_z[T]$ has distribution $\Theta_{\pi^b,z}(\mu^a_z)$
(resp. $\Theta_{\pi^b,z}(\mu^b_z)$),
where $\mu^a_z$ is solution to the RDE:
$$\mu^a_z=\Theta_{\widehat{\pi}^a,z}\circ
\Theta_{\widehat{\pi}^b,z}(\mu^a_z),$$
and $\mu^b_z=\Theta_{\widehat{\pi}^b,z}(\mu^a_z)$.

For $z=0$: conditionally on the root being of type $a$
(resp. $b$), the random variable $\cR_*[T]$ has law
$\Theta_{\pi^a,\widehat{\pi}^b}(\mu_*^a)$ (resp. $\Theta_{\pi^b,\widehat{\pi}^a}(\mu_*^b)$), where $\mu^a_*$ is the largest solution to the RDE
\begin{equation}
\label{eq:rdezeroa}
\mu^a_*=\Theta_{\widehat{\pi}^a,\widehat{\pi}^b}(\mu^a_*),
\end{equation}
and $\mu^b_*$ is the largest solution to the RDE
$\mu^b_*=\Theta_{\widehat{\pi}^b,\widehat{\pi}^a}(\mu^b_*)$.
\end{lemma}

We now analyze the RDE (\ref{eq:rdezeroa}).
We define:
\begin{eqnarray}
\label{eq:defa}F^a(x) &=& \phi^a\left(
  1-\hphi^b(1-x)\right)-\frac{\phi^a{}'(1)}{\phi^b{}'(1)}\left(1-\phi^b(1-x)-x\phi^b{}'(1-x).
\right)
\end{eqnarray}
Observe that
\begin{eqnarray*}
F^a{}'(x) = \frac{\phi^a{}'(1)}{\phi^b{}'(1)}\phi^b{}''(1-x)\left(\hphi^a(1-\hphi^b(1-x)) -x\right).
\end{eqnarray*}
Hence any $x$ where $F^a$ admits a local maximum must satisfy
$x=\hphi^a(1-\hphi^b(1-x))$. We define the historical records of $F^a$
as the set of $x \in [0,1]$ such that $x=\hphi^a(1-\hphi^b(1-x))$ and  for any $0 \leq y < x$, $F^a (x) > F^a (y)$ (the latter condition being empty if $x = 0$).  

\begin{theorem}
\label{th:rdea}
If $p_1<\ldots<p_r$ are the locations of the historical records of $F^a$, then the RDE (\ref{eq:rdezeroa}) admits exactly $r$ solutions ; moreover, these solutions can be stochastically ordered, say $\mu_{1}{<}\ldots{<}\mu_{r}$, and for any $i\in\{1,\ldots,r\}$,
\begin{itemize}
\item $\mu_{i}((0,1])=p_i$ ;
\item $\Theta_{\pi^a,\widehat{\pi}^b}(\mu_{i})$ has mean $F^a(p_i)$.
\end{itemize}
\end{theorem}
The proof of Theorem \ref{th:rdea} relies on two lemmas.
\begin{lemma}
\label{lm:contincra}
The operators $\Theta_{\pi^a,\widehat{\pi}^b}$ and
$\Theta_{\widehat{\pi}^a,\widehat{\pi}^b}$ are continuous (with
respect to weak convergence) and strictly increasing (with respect to
stochastic ordering) on $\cP\left([0,1]\right)$.
\end{lemma}
\bp{ of Lemma \ref{lm:contincra}}
It follows directly from the fact that, for any $n\geq 0$ and any
$n_1,\ldots,n_n\geq 0$, the
mapping $$x\mapsto\frac{1}{1+\sum_{i=1}^n\left(\sum_{j=1}^{n_i}x_{ij}\right)^{-1}}$$
is continuous and increasing from $[0,1]^{n_1+\ldots+n_n}$ to $[0,1]$.
\ep
\begin{lemma}
\label{lm:pFa(p)}
For any $\mu\in\cP\left([0,1]\right)$, letting $p=\mu\left((0,1]\right)$, we have
\begin{enumerate}
\item $\Theta_{\widehat{\pi}^a,\widehat{\pi}^b}(\mu)\left((0,1]\right)=\hphi^a(1-\hphi^b(1-p))$
\item if $\Theta_{\widehat{\pi}^a,\widehat{\pi}^b}(\mu)\leq\mu$, then the mean of $\Theta_{{\pi}^a,\widehat{\pi}^b}(\mu)$ is at least $F^a(p).$
\item if $\Theta_{\widehat{\pi}^a,\widehat{\pi}^b}(\mu)\geq\mu$, then the mean of $\Theta_{{\pi}^a,\widehat{\pi}^b}(\mu)$ is at most $F^a(p)$;
\end{enumerate}
In particular, if $\mu$ is a fixed point of $\Theta_{\widehat{\pi}^a,\widehat{\pi}^b}$, then $p=\hphi^a(1-\hphi^b(1-p))$ and $\Theta_{{\pi}^a,\widehat{\pi}^b}(\mu)$ has mean $F^a(p)$.
\end{lemma}
\bp{ of Lemma \ref{lm:pFa(p)}}
In equation (\ref{eq:operator}) it is clear that $Y>0$ if and only if
for any $i\in\{1,\ldots,\cN\}$, there exists
$j\in\{1,\ldots,{\cN_i}'\}$ such that $X_{ij}>0$. With the notation
introduced above, this rewrites:
$$\Theta_{\widehat{\pi}^a,\widehat{\pi}^b}(\mu)\left( (0,1]\right)=\widehat{\phi}^a\left(1-\widehat{\phi}^b\left(1-\mu\left((0,1]\right)\right)\right),$$
hence the first result follows.

Now let $X\sim\mu$, $Y\sim
\Theta_{\widehat{\pi}^a,\widehat{\pi}^b}(\mu)$, $\cN^a\sim\pi^a$,
$\widehat{\cN}^a\sim\widehat{\pi}^a$, and let $S,S_1,\ldots$ have the
distribution of the sum of a $\widehat{\pi}^b-$distributed number of
i.i.d. copies of $X$, all these variables being independent. Observe that
\begin{eqnarray*}
\frac{1}{1+\sum_{i=1}^{\cN^a} {S_i}^{-1}} 
& = & \left(1-\frac{\sum_{i=1}^{\cN^a} {S_i}^{-1}}{1+\sum_{i=1}^{\cN^a} {S_i}^{-1}}\right) \ind_{\left\{\forall i=1\ldots\cN^a, S_i>0\right\}} \\
& = &  \ind_{\left\{\forall i=1\ldots\cN^a, S_i>0\right\}} -\sum_{j=1}^{\cN^a} \frac{ {S_j}^{-1}}{1+S_{j} ^{-1}  + \sum_{ 1 \leq i \leq \cN^a , i \ne j} {S_i}^{-1}} \ind_{ \left\{\forall i=1\ldots\cN^a, S_i>0\right\}}
\end{eqnarray*}
Then, $\Theta_{\pi^a,\widehat{\pi}^b}(\mu)$ has mean
\begin{eqnarray*}
\EE\left[\frac{1}{1+\sum_{i=1}^{\cN^a} {S_i}^{-1}}\right] & = & \PP \left( \forall i=1\ldots\cN^a, S_i>0\right) - \sum_{k = 1 } ^\infty k \pi^a_k \EE\left[\frac {S^{-1}}{S^{-1}+1+\sum_{i=1}^{k-1}{S_i}^{-1}}\mathbf 1_{\left\{S>0,\forall i=1\ldots k -1, S_i>0\right\}}\right] \\
& = & \phi^a(1-\widehat{\phi}^b(1-p))-\phi^a{}'(1)\EE\left[\frac {S^{-1}}{S^{-1}+1+\sum_{i=1}^{\widehat{\cN}^a}{S_i}^{-1}}\mathbf 1_{\left\{S>0,\forall i=1\ldots\widehat{\cN}^a, S_i>0\right\}}\right]\\
& = &
\phi^a(1-\widehat{\phi}^b(1-p))-\phi^a{}'(1)\EE\left[\frac{Y}{Y+S}\mathbf 1_{\left\{S>0\right\}}\right],
\end{eqnarray*}
where the second and last lines follow from (\ref{eq:F}) and
$Y\sim\Theta_{\widehat{\pi}^a,\widehat{\pi}^b}(\mu)$, respectively.
Now, for any $s>0$, $x\mapsto {x}/{x+s}$ is increasing and hence, depending on whether $\Theta_{\widehat{\pi},\widehat{\pi}}(\mu)\geq\mu$ or $\Theta_{\widehat{\pi}^a,\widehat{\pi}^b}(\mu)\leq\mu$, $\Theta_{\pi^a,\widehat{\pi}^b}(\mu)$ has mean at most/least:
\begin{eqnarray*}
\label{eq:calc}
\phi^a(1-\widehat{\phi}^b(1-p))-\phi^a{}'(1)\EE\left[\frac{X}{X+S}\mathbf
  1_{\left\{S>0\right\}}\right]   =
\phi^a(1-\widehat{\phi}^b(1-p))-\phi^a{}'(1)\EE\left[\frac{X  }{X  +\sum_{ i = 1} ^{\widehat{\cN}^b} X_i  }\ind_{\left\{\cN^*\geq 1\right\}} \right],
\end{eqnarray*}
with $X_i$ are i.i.d. copies of $X$ independent of $\widehat{\cN}^b\sim\widehat{\pi}^b$ and $\cN^*=\sum_{i=1}^{\widehat{\cN}^b}\mathbf 1_{\{X_i>0\}}$. Now if $X'$ is the law of $X$ conditioned on $\{X > 0\}$, and $X'_i$ are i.i.d. copies of $X'$, by exchangeability, we find
$$
\EE\left[\frac{X  }{X +\sum_{ i = 1} ^{\widehat{\cN}^b} X_i  }\ind_{\left\{\cN^*\geq 1\right\}} \right] = p \EE\left[\frac{X'  }{X'+\sum_{ i = 1} ^{ \cN^*} X'_i  }\ind_{\left\{\cN^*\geq 1\right\}} \right]  = p \EE\left[\frac{1}{1+\cN^*}\mathbf 1_{\left\{\cN^*\geq 1\right\}}\right].
$$
Hence finally, depending on whether $\Theta_{\widehat{\pi},\widehat{\pi}}(\mu)\geq\mu$ or $\Theta_{\widehat{\pi}^a,\widehat{\pi}^b}(\mu)\leq\mu$, $\Theta_{\pi^a,\widehat{\pi}^b}(\mu)$, $\Theta_{\pi^a,\widehat{\pi}^b}(\mu)$ has mean at most/least:
\begin{eqnarray*}
\phi^a(1-\widehat{\phi}^b(1-p))-p\phi^a{}'(1)\EE\left[\frac{1}{1+\cN^*}\mathbf 1_{\left\{\cN^*\geq 1\right\}}\right]
\end{eqnarray*}
But using the definition (\ref{eq:F}) and the combinatorial identity $(n+1){n\choose d}=(d+1){n+1\choose d+1}$, one easily derive :
\begin{eqnarray*}
&&\phi^a(1-\widehat{\phi}^b(1-p))-p\phi^a{}'(1)\EE\left[\frac{1}{1+\cN^*}\mathbf 1_{\left\{\cN^*\geq 1\right\}}\right]\\
&=& \phi^a(1-\widehat{\phi}^b(1-p))- p\phi^a{}'(1)\sum_{n\geq 1}\widehat{\pi_n}^b\sum_{d=1}^n{n\choose d}\frac {p^d(1-p)^{n-d}}{d+1}
 =  F^a(p).
\end{eqnarray*}
\ep

\bp{ of Theorem \ref{th:rdea}}
Let $p\in[0,1]$ such that $\hphi^a(1-\hphi^b(1-p))=p$, and define
$\mu_0=\textrm{Bernoulli}(p)$.
From Lemma \ref{lm:pFa(p)} we know that
$\Theta_{\widehat{\pi}^a,\widehat{\pi}^b}\left(\mu_0\right)\left((0,1]\right)=p$,
and since \textrm{Bernoulli(p)} is the largest element of $\cP([0,1])$
putting mass $p$ on $(0,1]$, we have
$\Theta_{\widehat\pi^a,\widehat\pi^b}\left(\mu_0\right)\leq\mu_0$.
Immediately, Lemma \ref{lm:contincra} guarantees that the limit
$$\mu_\infty=\lim_{k\to\infty}\searrow\Theta^k_{\widehat\pi^a,\widehat\pi^b}\left(\mu_0\right)$$
exists in $\cP\left([0,1]\right)$ and is a fixed point of
$\Theta_{\widehat\pi^a,\widehat\pi^b}$. Moreover, by Fatou's lemma,
the number $p_\infty=\mu_\infty\left((0,1]\right)$ must satisfy
$p_\infty\leq p$. But then the mean of
$\Theta_{\pi^a,\widehat\pi^b}(\mu_\infty)$ must be both
\begin{itemize}
\item equal to $F^a(p_\infty)$ by Lemma \ref{lm:pFa(p)} with $\mu=\mu_\infty$ ;
\item at least $F^a(p)$ since this holds for all
  $\Theta_{\pi^a,\widehat\pi^b}\circ
  \Theta^k_{\widehat\pi^a,\widehat\pi^b}(\mu_0)$, $k\geq 1$ (Lemma \ref{lm:pFa(p)} with $\mu=\Theta^k_{\widehat\pi^a,\widehat\pi^b}(\mu_0)$).
\end{itemize}
We have just shown both $F^a(p)\leq F^a(p_\infty)$ and $p_\infty\leq
p$. From this, we now deduce the one-to-one correspondence between
historical records of $F^a$ and fixed points of
$\Theta_{\widehat\pi^a,\widehat\pi^b}$.
We treat each inclusion separately:
\begin{enumerate}
\item If $F^a$ admits an historical record at $p$, then clearly $p_\infty=p$, so $\mu_\infty$ is a fixed point satisfying $\mu_\infty\left((0,1]\right)=p$.
\item Conversely, considering a fixed point $\mu$ with
  $\mu\left((0,1]\right)=p$, we want to deduce that $F^a$ admits an
  historical record at $p$. We first claim that $\mu$ is the above
  defined limit $\mu_\infty$. Indeed, $\mu\leq\textrm{Bernoulli(p)}$
  implies $\mu\leq\mu_\infty$ ($\Theta_{\widehat\pi^a,\widehat\pi^b}$
  is increasing), and in particular $p\leq p_\infty$. Therefore,
  $p=p_\infty$ and $F^a(p)=F^a(p_\infty)$. In other words, the two
  ordered distributions
  $\Theta_{\pi^a,\widehat\pi^b}(\mu)\leq\Theta_{\pi^a,\widehat\pi^b}(\mu_\infty)$
  share the same mean, hence are equal. This ensures
  $\mu=\mu_\infty$.
Now, if $q<p$ is  any historical record location, we know from part 1 that
$$\nu_\infty=\lim_{k\to\infty}\searrow
\Theta_{\widehat\pi^a,\widehat\pi^b}^{k}\left(\textrm{Bernoulli(q)}\right)$$
is a fixed point of $\Theta_{\widehat\pi^a,\widehat\pi^b}$ satisfying
$\nu_\infty\left((0,1]\right)=q$. But $q<p$, so
$\textrm{Bernoulli(q)}<\textrm{Bernoulli(p)}$, hence
$\nu_\infty\leq\mu_\infty$. Moreover, this limit inequality is strict
because $\nu_\infty\left((0,1]\right) = q < p =
\mu_\infty\left((0,1]\right)$. Consequently,
$\Theta_{\pi^a,\widehat\pi^b}(\nu_\infty)<
\Theta_{\pi^a,\widehat\pi^b}(\mu_\infty)$ and taking expectations, $F^a(q)<
F^a(p)$. Thus, $F^a$ admits an historical record at $p$.
\end{enumerate}
\ep

We may now finish the proof of Theorem \ref{th:KSbi}.

\noindent{\em Proof of Theorem \ref{th:KSbi} : case of bounded degrees. }
We assume that $\pi_a$ and $\pi_b$ have bounded support. Recall
(\ref{eq:contrlamb}), so that we have
$\lambda=\frac{\phi^b{}'(1)}{\phi^a{}'(1)+\phi^b{}'(1)}$, where $\lambda$ is the probability that the root is of type $a$.
Theorems \ref{co:cvMN} and \ref{th:rdea} and Lemma \ref{lem:RDEbi} give:
\begin{equation}
\label{eq:KSbi2}
\frac{\nu(G_n)}{{|V_n|}}\xrightarrow[n\to\infty]{}
\frac{\lambda(1-\max_{x\in[0,1]}F^a(x))+(1-\lambda)(1-\max_{x\in[0,1]}F^b(x))}{2},
\end{equation}
where $F^a$ is defined in (\ref{eq:defa}) and $F^b$ is defined similarly by
 \begin{eqnarray}
\label{eq:defb}F^b(x) &=& \phi^b\left(1-\hphi^a(1-x)\right)-\frac{\phi^b{}'(1)}{\phi^a{}'(1)}\left(1-\phi^a(1-x)-x\phi^a{}'(1-x)
\right).
\end{eqnarray}
For any $x$ which is an historical record of $F^a$, we define
$y=\hphi^b(1-x)$ so that $\hphi^a(1-y)=x$. Then we have:
\begin{eqnarray*}
\lambda(1-F^a(x)) &=& \lambda\left(1-\phi^a(1-y)+\phi^a{}'(1)\left(
    \frac{1}{\phi^b{}'(1)}-\frac{\phi^b(1-\hphi^a(1-y))}{\phi^b{}'(1)}-y\hphi^a(1-y)\right)\right)\\
&=& (1-\lambda)(1-F^b(y)).
\end{eqnarray*}
By symmetry, this directly
implies that
$\lambda(1-\max_{x\in[0,1]}F^a(x))=(1-\lambda)(1-\max_{x\in[0,1]}F^b(x))$
so that (\ref{eq:KSbi2}) is equivalent to (\ref{eq:KSbi}). This proves Theorems  \ref{th:KS} and \ref{th:KSbi} for distributions with bounded support. \ep

\noindent{\em Proof of Theorem \ref{th:KSbi} : general case. } To keep notation simple, we only prove Theorem \ref{th:KS}. The following proof clearly extends to the case of UHGW trees. Let $G_1,G_2,\ldots$ be finite random graphs whose local weak limit is a
Galton-Watson tree $T$, and assume that the degree distribution $\pi$ of
$T$ (with generating function $\phi$) has a finite mean : $\phi'(1)=\sum_{n}n\pi_n<\infty$. 
For any rooted graph $G$ and any fixed integer $d\geq 1$, recall that $G^d$ is the graph obtained from $G$ by deleting all edges adjacent to a vertex $v$ whenever $\deg(v)>d$. Hence $T^d$ is a Galton-Watson tree whose degree distribution $\pi^{d}$ is defined by
$$\forall i\geq 0, \pi^{d}_i=\pi_i{\bf
    1}_{i\leq d}  + \ind_{i = 0} \sum_{k \geq d+1} \pi_k.$$
    By Theorem \ref{th:main}, Equation \eqref{eq:toshow} and our weaker version of
Theorem \ref{th:KS} for distributions with bounded support,
\begin{equation}
\label{eq:arg1}
\frac{\nu(G_n)}{|V_n|}\xrightarrow[n\to\infty]{} \lim_{d \to \infty} \min_{x\in[0,1]}g^d(x),
\end{equation}
with $\phi_d (x)=\sum_{k=0}^{d}\pi_kx^k$
and $$g^d(x)=1-\frac{1}{2}(1-x)\phi_d'(x)-\frac{1}{2}\phi_d(x)-\frac{1}{2}\phi_d\left(1-\frac{\phi_d'(x)}{\phi_d'(1)}\right).$$
Also, as $d\to\infty$, we have $\phi_d\to\phi$ and
$\phi'_d\to\phi'$ uniformly on $[0,1]$, so 
\begin{equation}
\label{eq:arg3}
\min_{x\in[0,1]}g^d(x)\xrightarrow[n\to\infty]{}\min_{x\in[0,1]}g(x),
\end{equation}
with $g(x)=1-\frac{1}{2}(1-x)\phi'(x)-\frac{1}{2}\phi(x)-\frac{1}{2}\phi\left(1-\frac{\phi'(x)}{\phi'(1)}\right).$
Finally, combining (\ref{eq:arg1}) and (\ref{eq:arg3}), we
easily obtain the desired 
\begin{equation*}
\label{eq:arg4}
\frac{\nu(G_n)}{|V_n|}\xrightarrow[n\to\infty]{}\min_{x\in[0,1]}g(x).
\end{equation*}

\ep

\subsection{Proof of Corollary \ref{cor:fr}}\label{sec:corfr}

Note that in Corollary \ref{cor:fr}, we divide $\nu(G_n)$ by
$|V^a_n|=\lfloor \alpha m\rfloor$ instead of $|V^a_n|+|V^b_n|=\lfloor \alpha m\rfloor+m$, so
that by Theorem \ref{th:KSbi}, we have
$\frac{\nu(G_n)}{|V^a_n|} \xrightarrow[n\to\infty]{}\min_{t\in [0,1]} 1-F^a(t)$.
We have $\phi^a(x)=x^k$, $\phi^b(x)=e^{\alpha k
  (x-1)}$ so that we have:
\begin{eqnarray*}
F^a(x) &=& \left( 1-e^{-k\alpha
    x}\right)^k-\frac{1}{\alpha}\left(1-e^{-k\alpha x}-k\alpha
  xe^{-k\alpha x} \right)\\
F^a{}'(x)&=& k^2\alpha e^{-k\alpha x}\left( \left(1-e^{-k\alpha
      x}\right)^{k-1}-x\right).
\end{eqnarray*}
Let $x^*$ be defined as in Corollary \ref{cor:fr} as the largest
solution to $x=\left(1-e^{-k\alpha x}\right)^{k-1}$.
It is easy to check (see Section 6 in \cite{lel12} for a more general
analysis) that
\begin{eqnarray*}
\min_{t\in [0,1]} 1-F^a(t)=\min\{1,1-F^a(x^*)\}.
\end{eqnarray*}
Setting $\xi^*=k\alpha x^*$, we have
$\frac{\xi^*}{k\alpha}=(1-e^{-\xi^*})^{k-1}$, so that
\begin{eqnarray*}
\min_{t\in [0,1]}
1-F^a(t)=\min\left\{1,1-\frac{1}{\alpha}\left(e^{-\xi^*}+\xi^*
    e^{-\xi^*}+\frac{\xi^*}{k}(1-e^{-\xi^*})-1\right)\right\}.
\end{eqnarray*}
Since $z\mapsto \frac{z(1-e^{-z})}{1-e^{-z}-z e^{-z}}$ is
increasing in $z$, we see that $\xi^*\geq \xi$ if and only if
$\alpha\geq \alpha_c$ and we get
\begin{eqnarray*}
\min_{t\in [0,1]}
1-F^a(t)=1-\ind(\alpha\geq \alpha_c)\frac{1}{\alpha}\left(e^{-\xi^*}+\xi^*
    e^{-\xi^*}+\frac{\xi^*}{k}(1-e^{-\xi^*})-1\right).
\end{eqnarray*}

\section*{Acknowledgement}
We would like to thank Andrea Montanari and Guilhem Semerjian for explaining us a
key idea for the proof of Lemma \ref{lm:unifctrl}, as well as Nikolaos
Fountoulakis, David Gamarnik, James Martin and Johan W\"astlund for
interesting discussions.

The authors acknowledge the support of the French Agence Nationale de la Recherche (ANR) under reference ANR-11-JS02-005-01 (GAP project)

\bibliographystyle{abbrv}
\bibliography{mat}

\section*{Appendix : uniqueness and non-uniqueness at zero temperature}

\begin{figure}[htb]
\begin{center}
 \includegraphics[angle=0, height=3cm]{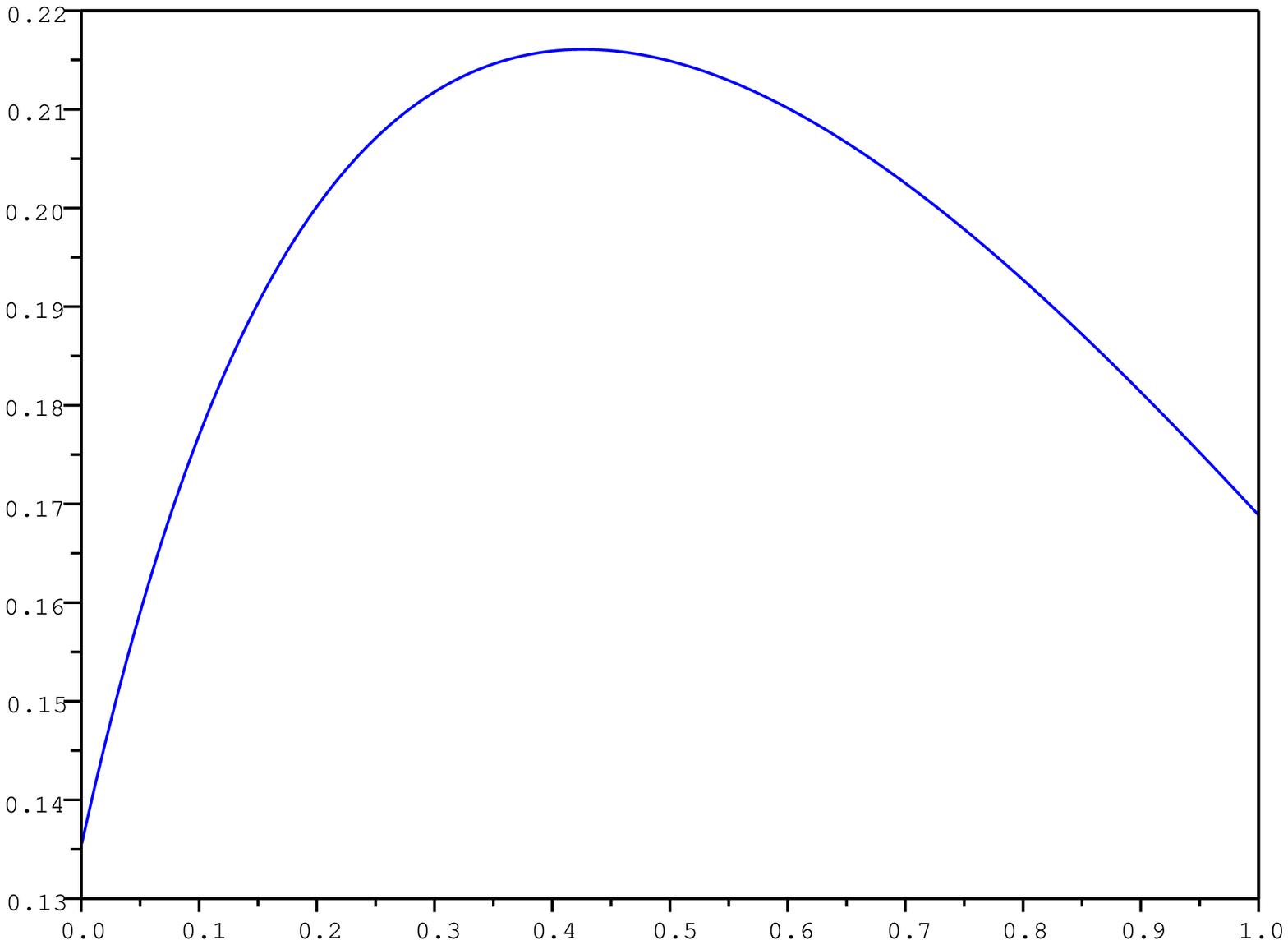} \hspace{30pt}
 \includegraphics[angle=0, height=3cm]{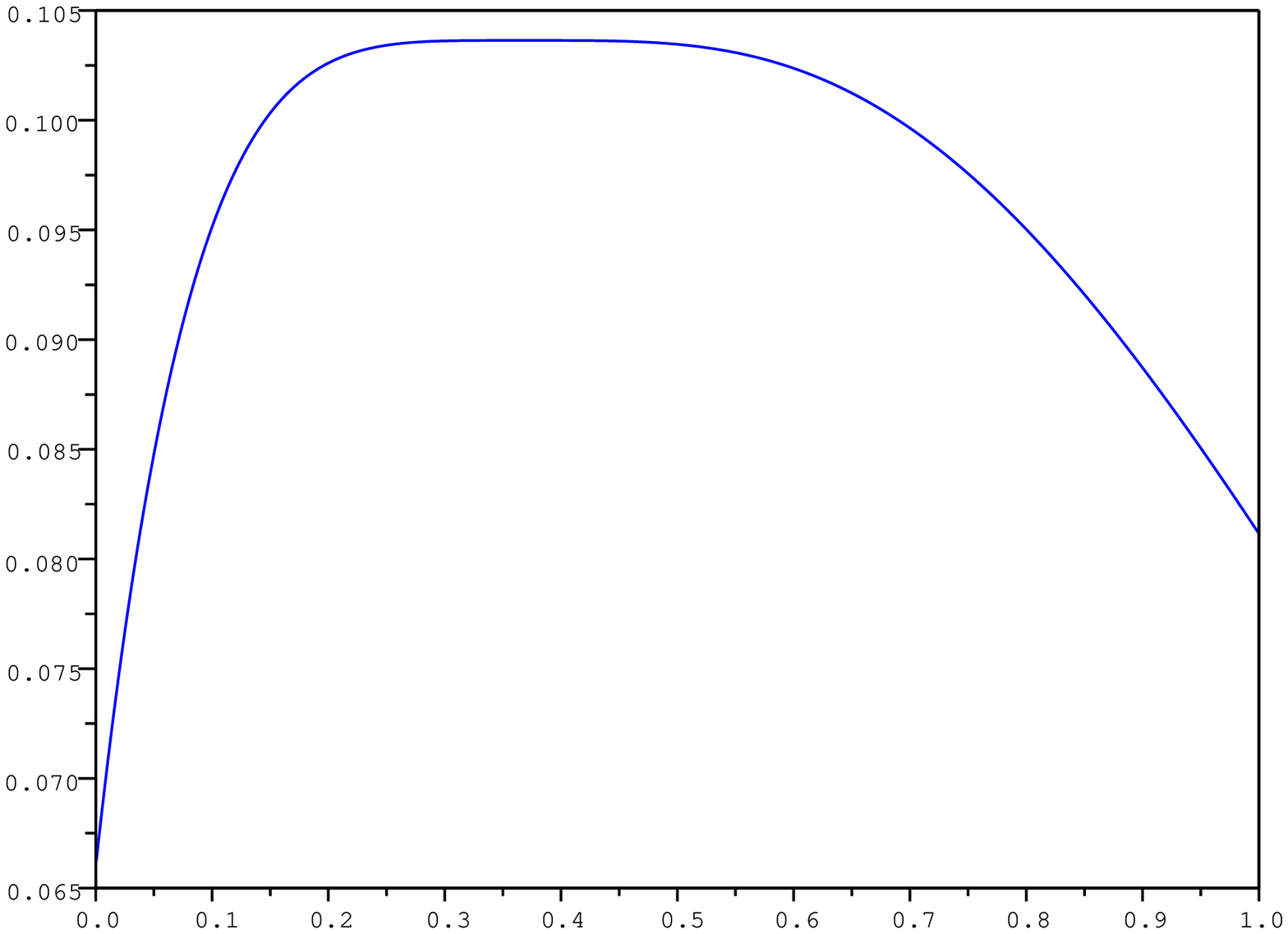} \hspace{30pt}
 \includegraphics[angle=0, height=3cm]{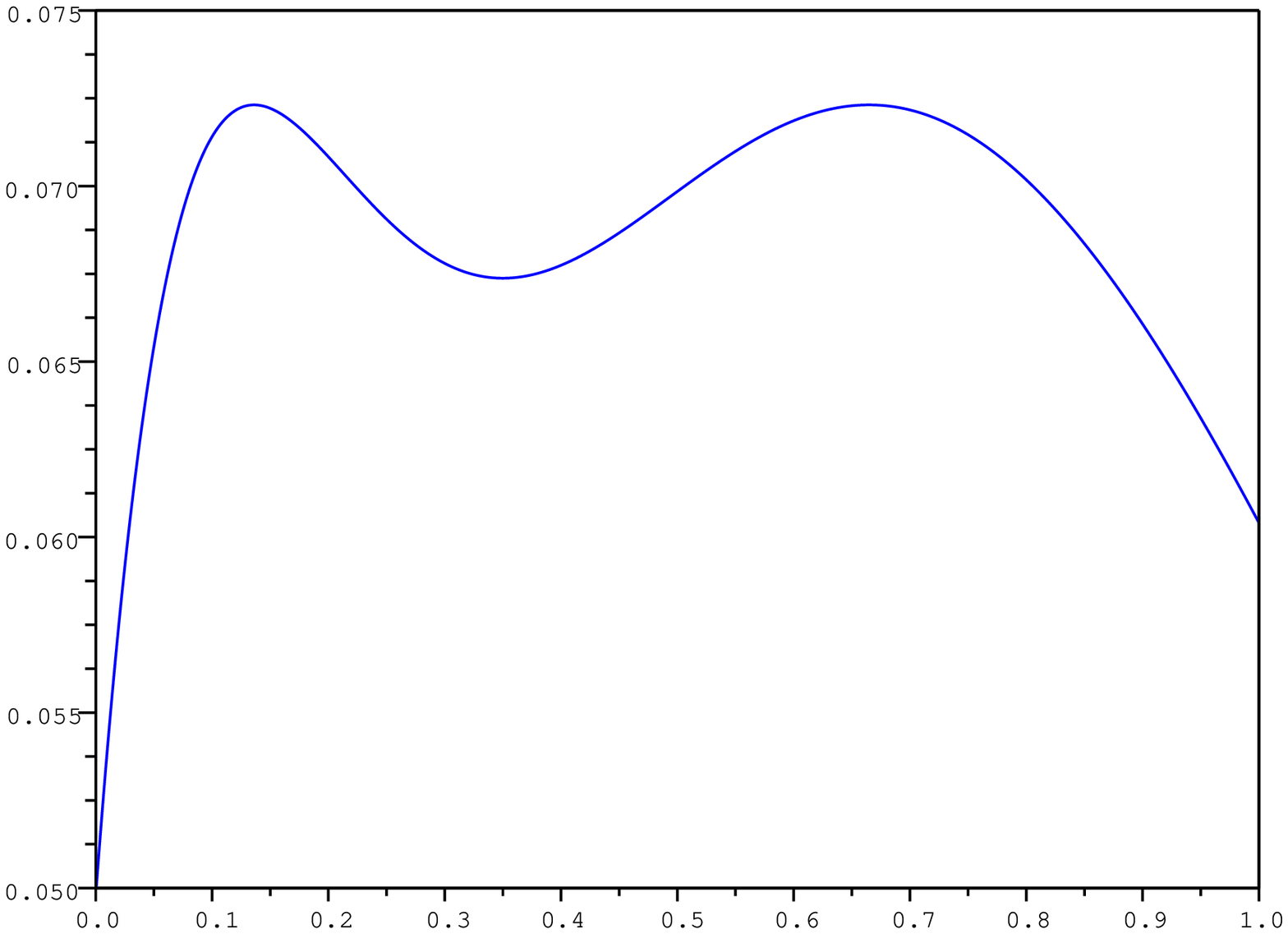} \hspace{30pt}
\caption{From left to right :  plot of $F$ for $c = 2$, $c = e$ and $c  = 3$. } \label{fig1}
\end{center}
\end{figure}

Consider a sequence $(G_n=(V_n,E_n))_{n\in \N}$ of finite graphs whose local weak
limit under uniform rooting is a UGW tree $T$. Let $\phi(t)=\sum\pi_{n}t^n$ be the generating function of the degree distribution $\pi$ of $T$, and for $t\in[0,1]$ set
\begin{eqnarray*}
F(t)=t\phi'(1-t)+\phi(1-t)+\phi\left(1-\frac{\phi'(1-t)}{\phi'(1)}\right)-1. 
\end{eqnarray*}
From Theorem \ref{th:rdea}, we know that there is a.s. a unique solution to the local recursion at temperature $z=0$ on $T$ (correlation decay) if and only if the first local extremum of $F$ is a global maximum. In that case, the convergence 
\begin{eqnarray}
\label{eq:finalcv}
\frac{\nu(G_n)}{|V_n|} \xrightarrow[n\to\infty]{}\frac{1-\max_{t\in [0,1]}F(t)}2
\end{eqnarray}
can be obtained by a fairly standard compactness-uniqueness argument, without any need for a detour through the positive temperature regime. As it is not hard to check, a sufficient condition for the first local extremum of $F$ to be a global maximum is that $\phi''$ is log-concave. This is in particular true in the Erd\H{o}s-R\'enyi case, where $\phi(t)=\exp(ct-1)$, $(c>0)$. The corresponding
function $F$ is given in Figure \ref{fig1} for various values of
$c$. 

However, there are simple examples of degree distributions $\pi$ for which the function $F$ has more than one historical record, implying the coexistence of multiple solutions to the local recursion (lack of correlation decay). In that case, a detour through the positive temperature regime is needed in order to establish the convergence (\ref{eq:finalcv}). Here is a couple of examples.
\begin{figure}[htb!]
\begin{center}
\includegraphics[angle=0, height=3cm]{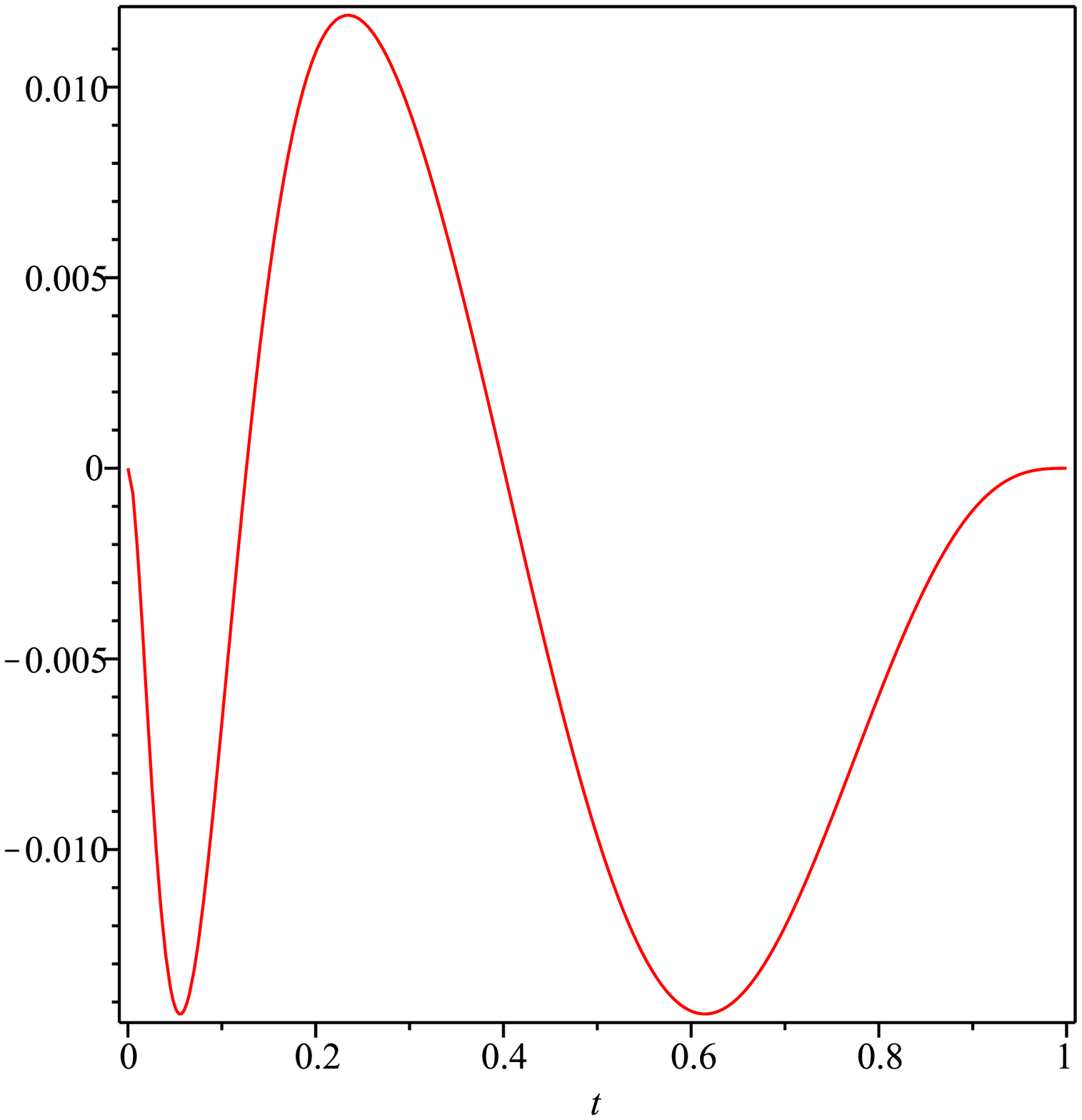} 
\hspace{30pt}
\includegraphics[angle=0, height=3cm]{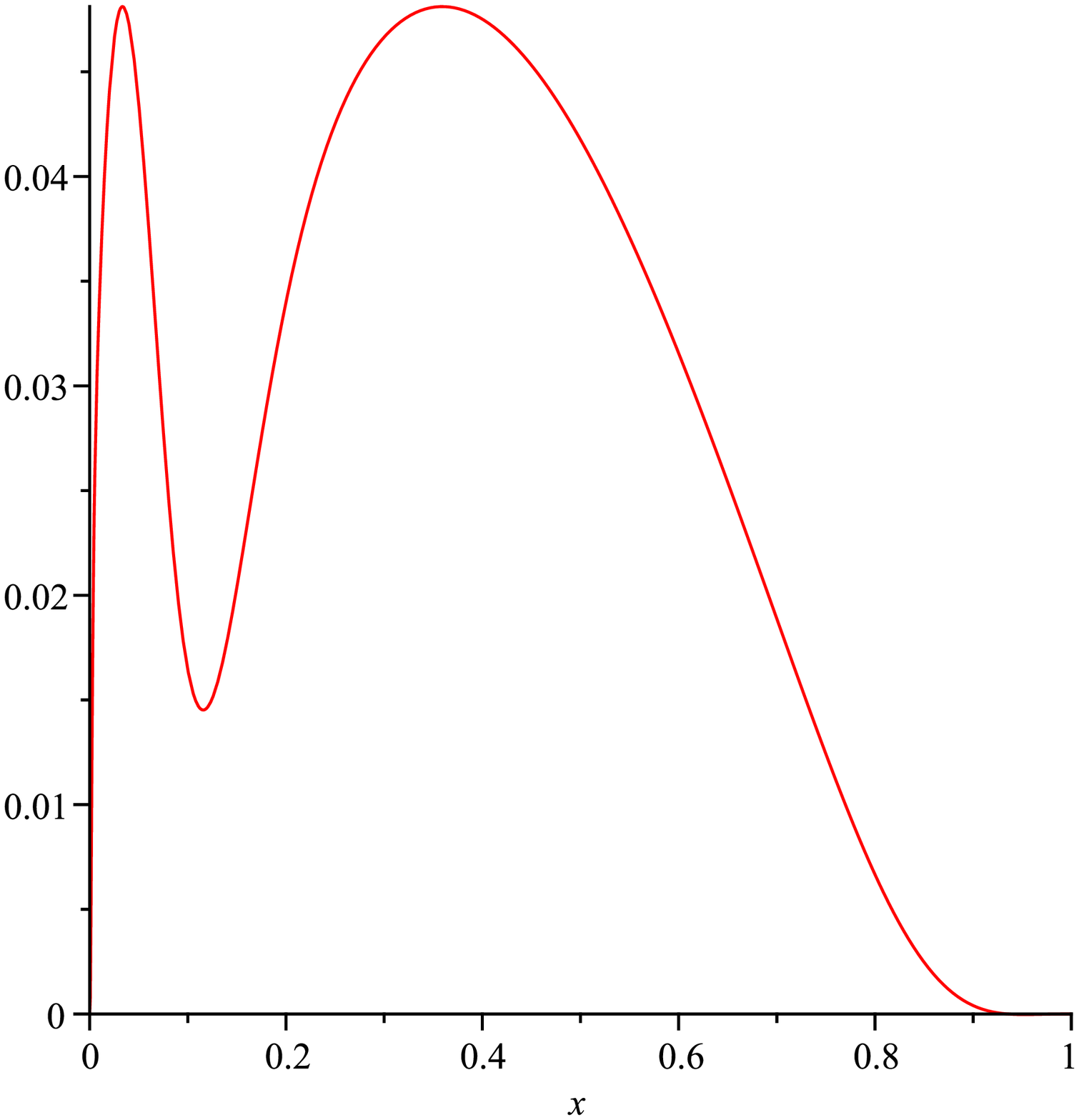}
\caption{Plots of $F$ for $\phi(t)=\frac{3}{4}t^3+\frac{1}{4}t^{15}$ (left) and $\phi(t)=\frac{50}{101}t^{3}+\frac{50}{101}t^{20}+\frac{1}{101}t^{700}$ (right).}\label{fig}
\end{center}
\end{figure}

\end{document}